\documentclass{article}

\usepackage{PRIMEarxiv}

\usepackage[utf8]{inputenc} % allow utf-8 input
\usepackage[T1]{fontenc}    % use 8-bit T1 fonts
\usepackage{hyperref}       % hyperlinks
\usepackage{url}            % simple URL typesetting
\usepackage{booktabs}       % professional-quality tables
\usepackage{amsfonts}       % blackboard math symbols
\usepackage{nicefrac}       % compact symbols for 1/2, etc.
\usepackage{microtype}      % microtypography
\usepackage{lipsum}
\usepackage{fancyhdr}       % header
\usepackage{graphicx}       % graphics
\graphicspath{{media/}}     % organize your images and other figures under media/ folder

\usepackage{amsfonts,amsmath,amssymb}
\usepackage{amsthm}
\newtheorem{theorem}{Theorem}[section]
\newtheorem{corollary}[theorem]{Corollary}
\newtheorem{lemma}[theorem]{Lemma}
\newtheorem{proposition}[theorem]{Proposition}

\theoremstyle{definition}
\newtheorem{definition}[theorem]{Definition}
\newtheorem{example}[theorem]{Example}
\newtheorem{remark}[theorem]{Remark}
\DeclareMathOperator{\Son}{Son} 
\DeclareMathOperator{\length}{length} 
\DeclareMathOperator{\card}{card}
\DeclareMathOperator{\Inv}{Inv}     % Invariant subspace
\DeclareMathOperator{\Hinv}{Hinv}

\newcommand{\CC}{\mathbb C}
\newcommand{\NN}{\mathbb N}

\title{Isomorphisms between lattices of hyperinvariant subspaces
}

\author{David Mingueza \\
  Nestlé Spain \\
  Esplugas de Llobregat\\
  \texttt{david.mingueza@outlook.es} \\
   \And
  M.Eulàlia Montoro\\
  Departament de Matemàtiques i Informàtica \\
  Universitat de Barcelona \\
  Barcelona\\
  \texttt{eula.montoro@ub.edu} \\
  \And
  Alicia Roca\\
  Departamento de Matemática Aplicada, IMM\\
  Universitat Politècnica de València\\
  València\\
  \texttt{aroca@mat.upv.es} \\
}  
\begin{document}
\maketitle
\begin{abstract} Given two nilpotent  endomorphisms,
we  determine when their lattices of hyperinvariant subspaces are isomorphic. The study of the lattice of hyperinvariant subspaces can be reduced to the nilpotent case   when the endomorphism has  a Jordan-Chevalley decomposition; for example,  it occurs if the underlying field is the field of complex numbers.
\end{abstract}

% keywords can be removed

%\fntext[fn1]{Partially supported by grant PID2019-104047GB-I00 and by grant PID2021-124827NB-I00 funded by MCIN/AEI/ 10.13039/501100011033 and by ``ERDF A way of making Europe'' by the ``European Union'' .}
%\fntext[fn3]{Partially supported by grant PID2021-124827NB-I00 funded by MCIN/AEI/ 10.13039/501100011033 and by ``ERDF A way of making Europe'' by the ``European Union'' .  }

\section{Introduction}~\label{sec:Intro}
Three main lattices of subspaces are naturally associated to endomorphisms, the lattices of invariant subspaces, of characteristic subspaces and of  hyperinvariant   subspaces. 
A vector subspace $V\subseteq \mathbb{F}^{n}$ is 
invariant with respect to the endomorphism $f$ of $\mathbb{F}^{n}$ if $f(V)\subset V$, the subspace is characteristic if it is invariant  for every automorphism commuting with $f$, and it is hyperinvariant if it is invariant for every  endomorphism commuting with $f$.

It is well-known when two lattices of invariant subspaces are isomorphic (\cite{BrickFill65}), and  isomorphisms of lattices of characteristic subspaces is still an open subject of study.  In this paper we characterize when two lattices of  hyperinvariant subspaces are isomorphic.
This problem was analyzed in  P. Y. Wu, "Which Linear Transformations Have Isomorphic Hyperinvariant Subspace Lattices?" (\cite{Wu92}). When studying the paper we realized that, although the statement of the main theorem is correct, some of the previous lemmas, theorems, or its proofs,  were not accurate, therefore they  need some reformulation and, as a consequence,   the  inclusion of the analysis of some missed cases. Hence, this work can be regarded as a review of  \cite{Wu92}.
 
The paper is organized as follows. In Section 2 we state some basic definitions and properties of lattices, in particular of lattices of invariant subspaces, and  of  the algebra of the centralizer. 
In Section 3 we describe the structure and some specific properties of the lattice of hyperinvariant subspaces, wich will be used later on.  Also, some necessary conditions for two lattices to be isomorphic are recalled.
In Section 4,  some special chains of hyperinvariant subspaces are presented, which will be a key  tool to study the problem. 
In Section \ref{isomorphisms} the special chains are widely used,  among other properties,  to study all possible isomorphisms of hyperinvariant lattices, in Subsection \ref{segrerneqs}  when the lattices have different Segre characteristic length, and 
in Subsection \ref{segrer=s} when they have the  same Segre characteristic length. Finally, in Subsection \ref{main} we summarize the results in a theorem.

\section{Preliminaries}
In this section we recall some results related to endomorphisms on complex vector spaces. We assimilate an endomorphism with its matrix representation.

A  partially ordered set $L$ such that every pair of elements $x,y\in L$ has a join ($x \vee y$) and a meet ($x \wedge y$) is called a lattice. A mapping between two lattices $L_1$ and $L_2$, $f: L_1 \rightarrow L_2$, is a lattice homomorphism if it preserves joints and meets.  If $f$ is an isomorphism we say that $L_1$ and $L_2$ are isomorphic and it is denoted by $L_1 \simeq L_2$.
If a lattice $L$ is finite, we denote by $\card(L)$ the number of elements of $L$.

%\medskip

Given $a,b \in L$ with $a\geq b$, the set formed by $\{x\in L: b \leq x \leq a\}$ is a sublattice of $L$ and is called the quotient sublattice $a/b$.

\medskip

An equivalent relation ($\sim$) on a lattice $L$ is a {\em congruence} relation if given  $a,b\in L$,  then 
\begin{equation} \label{def_congruence} 
a\sim b 
%\text{ implies } 
 \Rightarrow 
\left\{\begin{array}{l}(a\vee c) \sim (b\vee c),  \\
(a\wedge c) \sim (b\wedge c),
\end{array}\right. \text{   for every } c\in L.   
\end{equation}
We denote by 
$[a]=\{x\in L: x\sim a\}$ the class of elements related to $a$, and by $L / \hspace{-0.11cm}\sim =\{[a]: a\in L\}$  the {\em quotient lattice} formed by all the classes in  $L$. 
Moreover,   
\begin{equation}\label{chain_un}
a\sim b  \
 \Rightarrow 
%{\color{red} \text{the elements of the  the quotient sublattice } } (a\vee b)/(a\wedge b) \sim a.
 \forall \  c \in (a\vee b)/(a\wedge b), \ c \sim a.
\end{equation}
 
Let $f: L \longrightarrow L_1$ be a lattice  homomorphism. The congruence  relation defined by $x\sim y$ if $f(x)=f(y)$ is called the {\em kernel} of $f$ ($\ker f$).
The following isomorphism theorem is satisfied.

\begin{theorem} \label{teor:iso}{\rm(\cite[Theorem 1.5]{Gr23})}
    Let $L$ and $L_1$ be lattices.  Let $f: L \longrightarrow L_1$ be an homomorphism onto. Then, $L/\ker f$ is isomorphic to $L_1$. 
\end{theorem}

Let  $\CC^{n\times n}$  be the algebra of $n\times n$ matrices with entries in the complex field $\mathbb{C}$. 
%The algebra of polynomials in $A\in\CC^{n\times n}$ is denoted as 
% $$\mathbb{C}[A]=\{\sum_{i=1}^{l}a_{i}A^{i},a_{i}\in\mathbb{C}\},$$
%and 
Given $A\in\mathbb{C}^{n\times n}$, the algebra of the matrices that commute with $A$, i.e.,  the centralizer of $A$, denoted $Z(A)$, is 
$$Z(A)=\{B\in \CC^{n\times n}:  AB=BA\},$$  
the lattice of invariant subspaces of $A$, denoted $\Inv(A)$, is
$$\Inv(A)=\{F\subset \CC^{n\times 1}: AF\subseteq F\},$$ 
and the lattice of hyperinvariant subspaces of $A$, denoted $\Hinv(A)$, is
$$\Hinv(A)=\{F\subset \CC^{n\times 1}: BF\subseteq F,\,\forall B\in Z(A)\}.$$

%$ {\color{red}We will also write ${\cal A}_1\simeq {\cal A}_2$ if ${\cal A}_1$ and ${\cal A}_2$ are isomorphic  algebras/ si se queda, definir isomorpf of algebras.}

As mentioned before, isomorphisms of lattices of invariant subspaces are characterized. The result is stated in the next theorem.
\begin{theorem}
{\rm(\cite{BrickFill65})}\label{}
 Let $A,B\in \CC^{n\times n}$. Then, the following statements are equivalent:
 \begin{enumerate}
     
\item $Inv(A)$ and $ Inv(B)$ are isomorphic.
 %\text{ if and only if } 
 \item 
 $A$ and  $B$ have the same Jordan structure. 
 \item 
 $Z(A)$ and $ Z(B)$ are isomorphic.
 \end{enumerate}
\end{theorem}

%\begin{theorem}~\cite{GLR86}
%Let $A,B\in \CC^{n\times n}$ be nilpotent matrices. Then,
% $$\mathbb{C}[A]\simeq \mathbb{C}[B]
 %\Leftrightarrow 
 %\text{ if and only if } m_{A}=m_{B},$$
% where $m_{A}$ denotes the minimal polynomial of $A$.
%\end{theorem}

The following theorem gives a sufficient condition for two lattices of hyperinvariant subspaces to be isomorphic. In the next section we will see that it is not necessary.

\begin{theorem}\label{teo:suf}{\rm (\cite{Long76})}
 Let $A, B\in \mathbb{C}^{n\times n}$. If $\Inv(A)\simeq \Inv(B)$, then $\Hinv(A)\simeq \Hinv(B)$.
 \end{theorem}

%\subsection{Quotient Lattice}\label{sec:Quo}

%The definition of a quotient lattice allows us to reduce the size of the original lattice and apply induction method to prove some results (see~Theorem~\ref{teore:iso2}). For this,  we need to recall two important definitions on lattices and a key isomorphism theorem.

%%%%%%%%%%%%%%%%%%%%%%%%%%%%%%%%%%%%%%%%%%%%%%%%%%%%%%%%%%%%%%%%%%%%%%%%%

\section{The lattice of hyperinvariant subspaces}

 The study of the lattice of hyperinvariant subspaces can be reduced to the nilpotent case   when the endomorphism has  a Jordan-Chevalley decomposition (see~\cite{FillHL77}). As the reduction applies on the complex field, throughout the paper we will assume that $A\in \CC^{n\times n}$ is  nilpotent.

Let $A\in \CC^{n\times n}$ be a nilpotent matrix and $\alpha=(\alpha_{1},\ldots,\alpha_{r})\in\mathbb{N}^{r}$ its Segre characteristic, where $\alpha_{1}\geq\ldots\geq\alpha_{r}>0$.

Let $u=(u_1, \ldots, u_r)\in \NN^r$  such that 
\begin{align} \label{hypertuples1} 
u_{1}\geq \ldots\geq u_{r}\geq 0, \\
\label{hypertuples2} 
\alpha_{1}-u_{1}\geq\ldots\geq\alpha_{r}-u_{r}\geq 0,
\end{align} and define
$$V(\alpha_{1},\ldots,\alpha_{r})=\{u\in\mathbb{N}^{r}: \text{ satisfying } (\ref{hypertuples1}) \text{ and } (\ref{hypertuples2}) \}.$$

\vspace{0.2cm}

 Given $u,v\in V(\alpha_{1},\ldots,\alpha_{r})$, we define the following partial order:
     $$u\subset v\quad \mbox{if}\quad u_{i}\leq v_{i},\  i=1,\ldots,r.$$
     
\begin{proposition}{\rm(\cite{FillHL77})}\label{propo:hyp}
     $V(\alpha_{1},\ldots,\alpha_{r}) $ with the partial order relation ``$\subset$'' is a lattice.
The meet and the join are defined, respectively, as
$$u\wedge v=(\min(u_{1},v_{1}),\ldots,\min(u_{r},v_{r})),$$
$$u \vee v=(\max(u_{1},v_{1}),\ldots,\max(u_{r},v_{r})).$$
\end{proposition}

The elements of   $V(\alpha_{1},\ldots,\alpha_{r})$ are called  {\em hypertuples}, and we refer to the lattice $V(\alpha_{1},\ldots,\alpha_{r})$ as a {\em hyperlattice}.

\begin{proposition}{\rm(\cite{FillHL77})}
Given $A\in\mathbb{C}^{n\times n}$ a nilpotent matrix with Segre Characteristic $(\alpha_{1},\ldots,\alpha_{r})$, there exists a lattice isomorphism
$$f:\Hinv(A)\longrightarrow V(\alpha_{1},\ldots,\alpha_{r}).$$

\end{proposition}

\vspace{0.2cm}

 As a consequence of this proposition,  we identify the hyperinvariant subspaces with their corresponding hypertuples. If there is no possible confusion, given a Segre characteristic $\alpha=(\alpha_{1},\ldots, \alpha_{r})$ we will also denote the corresponding hyperlattice  by $V(\alpha)$.

 Unlike the lattice of invariant subspaces $\Inv(A)$, for finite dimensional matrices the lattice $\Hinv(A)$ is always finite, and its cardinality is known; we recall it in the next theorem.
\medskip

\begin{theorem} {\rm(\cite{FillHL77})} The number of hyperinvariant subspaces of a nilpotent matrix $A\in \CC^{n\times n}$ with Segre characteristic $\alpha=(\alpha_1, ..., \alpha_r)$ is
$$(\alpha_1 - \alpha_2 + 1) \dots (\alpha_{r-1 }-\alpha_r + 1)(\alpha_r + 1).$$
\end{theorem}

From now on, we will only consider
the reduced case of Segre characteristics, i.e., we will assume that $\alpha_{1}>\ldots>\alpha_{r}$, because there exists a lattice isomorphism between the hyperlattices of the non reduced and the reduced cases (see~\cite{Wu92}).  

\medskip

\begin{remark}\label{rem21}
The condition of Theorem~\ref{teo:suf} is not necessary. For example, let $A, B\in\mathbb{C}^{n\times n}$ be nilpotent matrices with Segre characteristics  $\alpha=(l,l-1)$ and $\beta=(2l-1)$, $l>0$,  respectively.  As we can see in the following diagram $V(\alpha)\simeq V(\beta)$, 

\begin{equation*} \label{iso2-1}
\begin{array}{cccc}
      (l,l-1) &&& (2l-1)\\
   |&&&|\\
   (l-1,l-1)&&& (2l-2)\\
   |&&&|\\
   (l-1,l-2)&&& (2l-3)\\
   |&&&|\\
   \vdots &&& \vdots\\
   |&&&|\\
   (1,0)&&& (1)\\
   |&&&|\\
   (0,0)&&& (0)\\
  \end{array}
  \end{equation*}

\noindent
but the corresponding lattices of invariant subspaces are not isomorphic (there exist at last two different invariant subspaces of $A$ of dimension $1$, whilst there exists only one of $B$).
\end{remark}

 Our aim is to characterize when two lattices of hyperinvariant subspaces are isomorphic. The general scheme of this paper follows that of  \cite{Wu92}, but  as announced in the Introduction, some results there were not accurate, either incomplete or defective. We fix here those flaws either completing the results or providing new proofs. Some concepts must be introduced.

%\begin{definition}%\label{defi:dist}
  %Given $V(\alpha_{1},\ldots,\alpha_{r})$ and $u, v\in V(\alpha_{1},\ldots,\alpha_{r})$, the distance between $u=(u_{1},\ldots,u_{r})$ and $v=(v_{1},\ldots,v_{r})$ is defined as
  %$$d(u,v)=\sum_{i=1}^{r}\mid u_{i}-v_{i}\mid.$$
 %\end{definition}

\medskip

We start by  defining the notion of  ``son'' of an element of a hyperlattice.

\begin{definition}\label{defi:son}
 Given $u=(u_1, \ldots, u_r)\in V(\alpha_{1},\ldots,\alpha_{r})$, we say that 
 $v=(v_1, \ldots, v_r)$  is a {\em son} of $u$ if there exists $j\in \{1, \ldots, r\}$ such that $v_i=u_i$ for $i\in \{1, \ldots, r\}\setminus \{j\}$,  $v_j=u_j-1$, and $v\in V(\alpha_{1},\ldots,\alpha_{r})$.  

  %Given $u,v\in V(\alpha_{1},\ldots,\alpha_{r})$, we say that $v$ is a {\em son} of $u$ if $v\subset u$ and there no exists any $w\in V(\alpha_{1},\ldots,\alpha_{r})$ such that $v\subset w\subset u $. 
  \end{definition}
  
 \medskip
 
 We denote by $\Son(u)$ the set of  sons of $u$.  Notice that if $v$ is a {\em son} of $u$ then $v\subset u$ and there is no any $w\in V(\alpha_{1},\ldots,\alpha_{r})$ such that $v\subsetneq w\subsetneq u$. Whenever $v\in\Son(u)$ we will say that $u$ is a {\em father} of $v$, and if $v,w \in \Son(u)$ we will say that $v$ and $w$ are {\em brothers}.

\medskip

\begin{definition}\label{defi:chain}
Given $V(\alpha_{1},\ldots,\alpha_{r})$, a {\em chain} $C$ of length $t$ from $w_1$ to $w_t$ is a sequence of hypertuples $w_{1},\ldots,w_{t} \in V(\alpha_{1},\ldots,\alpha_{r})$ such that 
 $w_{i+1}\in \Son(w_{i}), \ i=1,\ldots,t-1.$
 We will write the chain $C$ as 
$w_{1}- \dots - w_{t}$.

 In some occasions we will reorder the chain in the reverse way; i.e., a chain  $w_{1}- \dots - w_{t}$ from $w_1$ to $w_t$ will satisfy that $w_{i}\in \Son(w_{i+1}), \ i=1,\ldots,t-1$.

%A sequence of hypertuples $\{w_1, \ldots, w_n\} \in V(\alpha_{1},\ldots,\alpha_{r})$ is a {\em chain} if $w_1\subset  \ldots \subset  w_n$ and $w_i \in \Son (w_{i+1})$ for $i=1, \ldots, n-1$. 

  \end{definition}

 The following properties are satisfied.

\begin{proposition} \label{prop:sons} Let $f: V(\alpha) \longrightarrow V(\beta)$ be a hyperlattice  isomorphism, and $u,v \in V(\alpha)$. Then,
\begin{enumerate}
\item If $v\subset u$, then $f(v)\subset f(u)$.
\item If  $v \in \Son(u)$, then  $f(v) \in \Son(f(u))$. As a consequence, if $C$ is a  chain in $V(\alpha)$ of length $t$,  its image $f(C)$ is a chain in $V(\beta)$ of length $t$.
\item If  $\Son(u)=\{v\}$, then  $\Son(f(u))=\{f(v)\}$. 
\end{enumerate}
\end{proposition}

\begin{proof}
\begin{enumerate}
\item Observe that $v \vee u=u$. As $f$ preserve joins, $f(v) \vee f(u)=f(u)$,  hence $f(v)\subset f(u)$.
\item If $v \in \Son(u)$, then $f(v)\subset f(u)$. Assume that $f(v)$ is not a son of $f(u)$. It means that there exists $w\in V(\beta)$ such that $f(v)\subsetneq w \subsetneq f(u)$. Then, $v\subsetneq f^{-1}(w) \subsetneq u$, which implies that $v\notin \Son(u)$.
\item We know that  $f(v) \in \Son(f(u))$. If there exists $w\in \Son(f(u))$, $w\neq f(v)$, then $f^{-1}(w)\in \Son(u)$ and $f^{-1}(w)\neq v$, which is a contradiction.

\end{enumerate}

\end{proof}

%\begin{lemma}\label{lem:dist}Let $f: V(\alpha_{1},\ldots,\alpha_{r})\longrightarrow V(\beta_{1},\ldots,\beta_{s}),$ be a lattice isomorphism. Then, for every chain $C$ of length $t$ its image $f(C)$ has length $t$. \end{lemma}
% \begin{proof} The result follows from the fact that a lattice isomorphism the image of a son must be a son. \end{proof}

%\begin{corollary}
 %$V(\alpha_{1},\ldots,\alpha_{r})\simeq V(\beta_{1},\ldots,\beta_{s})\Rightarrow \alpha_{1}+\ldots+\alpha_{r}=\beta_{1}+\ldots+\beta_{s}$
% If the lattices  $V(\alpha_{1},\ldots,\alpha_{r})$ and $ V(\beta_{1},\ldots,\beta_{s})$ are isomorphic, then %\Rightarrow $\alpha_{1}+\ldots+\alpha_{r}=\beta_{1}+\ldots+\beta_{s}$.
%\end{corollary}

\medskip

In the following lemma we characterize the sons of a given  hypertuple.
\begin{lemma}\label{lem:son}
 Let $V(\alpha_{1},\ldots,\alpha_{r})$ be a hyperlattice  and $u=(u_{1},\ldots,u_{r})\in V(\alpha)$. Then:
 \begin{enumerate}
  \item $(u_{1}-1,u_{2},\ldots,u_{r})\in \Son(u)$ if and only if $ u_{1}>u_{2}$.
  \item Given $i\in\{2,\ldots, r-1\}$, then $(u_{1},\ldots,u_{i}-1,\ldots,u_{r})\in \Son(u) $ if and only if $ u_{i}>u_{i+1}$ and $\alpha_{i-1}-u_{i-1}>\alpha_{i}-u_{i}$.
  \item $(u_{1},u_{2},\ldots,u_{r}-1)\in \Son(u) $ if and only if $ u_{r}\neq 0$ and $\alpha_{r-1}-u_{r-1}>\alpha_{r}-u_{r}$.
 \end{enumerate}
\end{lemma}
\begin{proof}
Let $u\in V(\alpha)$. Notice that subtracting one from a single  component of $u$, if the resulting tuple belongs to $V(\alpha)$, we obtain a son. Therefore,
 \begin{enumerate}
     \item $(u_{1}-1,u_{2},\ldots,u_{r})\in V(\alpha)$ if and only if $ u_{1}-1\geq u_{2}$ and $ \alpha_{1}-(u_{1}-1)\geq \alpha_{2}-u_{2}$.  As  $u\in V(\alpha)$, the  two conditions are equivalent to $u_{1}>u_{2}$.
  \item Given $i\in\{2,\ldots, r-1\}$, then $(u_{1},\ldots,u_{i}-1,\ldots,u_{r})\in V(\alpha) $  if and only if  $u_{i-1}\geq u_{i}-1\geq u_{i} $ and $\alpha_{i-1}-u_{i-1}\geq \alpha_{i}-(u_{i}-1)\geq \alpha_{i}-u_{i}$. As $u\in V(\alpha)$,  these  conditions are equivalent to $u_{i}>u_{i+1}$ and $\alpha_{i-1}-u_{i-1}>\alpha_{i}-u_{i}$.
  \item $(u_{1},u_{2},\ldots,u_{r}-1)\in \Son(u) $   if and only if  $ u_r>0,  u_{r-1}\geq u_{r}-1\quad \mbox{and}\quad  \alpha_{r-1}-u_{r-1}\geq \alpha_{r}-(u_{r}-1)$.
   As $u\in V(\alpha)$, these  conditions are equivalent to that 
  $u_{r}\neq 0$ and $\alpha_{r-1}-u_{r-1}>\alpha_{r}-u_{r}$.
 \end{enumerate}

\end{proof}

%\noindent
%{\color{blue}{Nota: Taking $\alpha_0=\alpha_1+1, u_0=u_1$ and $u_{r+1}=0$, the three cases reduce to the second one.}}

\begin{example}
  Let $u=(3,2,1), v=(4,2,1)\in V(5,3,1)$. Then,
  $$\Son(u)=\{(2,2,1),(3,1,1),(3,2,0)\} , \quad \Son(v)=\{(3,2,1), (4,2,0)\}.$$
 \end{example}

\begin{corollary}\label{cor:u=0}
 Let $u=(u_{1},\ldots,u_{r})\in V(\alpha)$. Then, 
 $$\Son(u)=\{\emptyset\} \  \Leftrightarrow \text{ if and only if  } \   u=(0,\ldots,0).$$
\end{corollary}

%\begin{corollary}
% Let $u=(u_{1},\ldots,u_{r})\in V(\alpha)$. Then,
 %$$ \card (\Son(u))= r \ \text{ if and only if } \
%\Leftrightarrow$$
 %$$ 0< u_{i-1}-u_{i}<\alpha_{i-1}-\alpha_{i},\quad i=2,\ldots, r.$$
%\end{corollary}

As a consequence of Corollary \ref{cor:u=0} and  Proposition \ref{prop:sons}, item 2, we obtain the following result.

\begin{corollary} 
\begin{enumerate}
\item Given $u=(u_1, \ldots, u_r) \in V(\alpha_1, \ldots, \alpha_r)$, there exists a chain $C$ 
from $u$ to  $0$ $(C: u- \dots - 0)$
of length $t=u_1+\dots+u_r+1$.
\item Let $f:V(\alpha_1, \ldots, \alpha_r) \longrightarrow V(\beta_1, \ldots, \beta_s)$ be a hyperlattice  isomorphism, and $v=f(u), u \in V(\alpha)$. Then,
$$u_1+\dots+u_r=v_1+\dots +v_s.$$
\end{enumerate}
\end{corollary}

%\begin{proof} Its suffices to consider a chain $C:u-\ldots-0$ and the fact that a son is obtained by substracting 1.
 %for any isomorphism of lattices of hypertuples, must be that 
 %\end{proof}

\begin{remark} \label{rem:sum-card}
 If the hyperlattices  $V(\alpha_{1},\ldots,\alpha_{r})$ and $ V(\beta_{1},\ldots,\beta_{s})$ are isomorphic, then the following two conditions must be satisfied:
 \begin{enumerate}
     \item \label{rem:sum}
%\begin{equation*} %\label{size}
Dimension condition: $\alpha_{1}+\dots+\alpha_{r}=\beta_{1}+\dots+\beta_{s}$.
%\end{equation*}

\item \label{rem:card}
%\begin{equation*} %\label{card}
%(\alpha_{1}-\alpha_{2}+1)\dots (\alpha_{r-1}-\alpha_{r-1}+1)(\alpha_r+1)=(\beta_{1}-\beta_{2}+1)\dots (\beta_{s-1}-\beta_{s-1}+1)(\beta_s+1).
Cardinality condition: $\card(V(\alpha_{1},\ldots,\alpha_{r}))=\card( V(\beta_{1},\ldots,\beta_{s}))$.
%\end{equation*}
\end{enumerate}

Along the paper we will frequently refer to these two conditions.
\end{remark}

 \section{Special Chains and Riding Special Chains}\label{sec:special}
 Useful tools  to study when two hyperlattices are isomorphic, are the notions of ``special chains'' and ``riding special chains''.
To define them we need to  characterize when a hypertuple has a unique son. 
 
 \noindent

\begin{theorem}\label{teo:son1}
 Let $u\in V(\alpha_{1},\ldots,\alpha_{r})$.  Let $k=\max\{i:u_1=\dots=u_i\}$  and $q=\max\{i: u_i> 0\}$.
 Then,
 $$\card(\Son(u))=1 \quad \text{if and only if} \quad
 %\Leftrightarrow  
 \alpha_{k}-u_{k}=\alpha_{k+1}-u_{k+1}=\ldots=\alpha_{q}-u_{q}.$$
  In such a case, $$\Son(u)=\{(u_{k},\ldots,u_{k}-1,u_{k+1},\ldots,u_{q},0,\ldots,0)\}.$$
\end{theorem}

\begin{proof}
 %Assume that $\#(\Son(u))=1$. 
If $k=r$ the result is trivial. If $k<r$, we have $u_{k}>u_{k+1}$, and 
$$\bar u=(u_{k},\ldots,u_{k}-1,u_{k+1},\ldots,u_{q},0\ldots,0)\in \Son(u).$$

Let $\card(\Son(u))=1$. If $k=q$, the result is trivial. If $k<q$, assume that there exists $k\leq j
 <q$ such that $\alpha_{j}-u_{j}>\alpha_{j+1}-u_{j+1}=\ldots=\alpha_{q}-u_{q}$. Then, by Lemma~\ref{lem:son}
$$\bar{\bar{u}}=(u_{k},\ldots,u_{k},\ldots, u_{j+1}-1,\ldots,u_{q},0,\ldots,0)\in \Son(u).$$
As $\bar{\bar{u}}\neq \bar u$, 
%$$(u_{k},\ldots,u_{k},\ldots, u_{j+1}-1,\ldots,u_{r}) \neq \bar u,$$
 $u$ has more than one son, which is a contradiction.

Conversely. Assume that $u$ has  a second son $v=(v_{1},\ldots,v_{r})$. It means that there exists only one index $1\leq j\leq r$ such that $v_j=u_{j}-1$. Then, $k\leq j\leq q$ ($j\geq k$, otherwise $v_j=u_j-1<u_{j+1}=v_{j+1}$. Also, $j\leq q$, otherwise $u_j=0$). If $k<j\leq q$, then, 
$$\alpha_{j-1}-u_{j-1}=\alpha_{j-1}-v_{j-1}\geq \alpha_{j}-v_{j}=\alpha_{j}-(u_{j}-1)>\alpha_{j}-u_{j},$$  
which by hypothesis cannot occur. Therefore, $j=k$ and $v=\bar u$.

\end{proof}

 \begin{example} \
 Let $k$ and $q$ be defined as in the previous theorem.
 \begin{enumerate} 
   \item Let $\alpha=(7,3,1)$. Then,
   $\Son(2,2,0)=\{(2,1,0)\}.$
  Here $k=q=2$.
 \item  Let $\alpha=(7,5,3,1)$. Then,
   $\Son(3,3,1,0)=\{(3,2,1,0)\}.$
   Here $k=2$ and $q=3$.
  
   \item Let $\alpha=(7,5,4,3)$. Then,
   $\Son(3,3,2,1)=\{(3,2,2,1)\}$.
   Here $k=2$ and $q=4$.
   \item Notice that  $\Son(\alpha_{1},\ldots,\alpha_{r})=\{(\alpha_{1}-1,\ldots,\alpha_{r})\}$.  Here $k=1$ and $q=r$.  
   \end{enumerate}
  \end{example}

%\begin{definition}
% Given $V(\alpha_{1},\ldots,\alpha_{r})$. A {\em chain} $C$ of length $t$ is a succession of hypertuples $w_{1},\ldots,w_{t} \in V(\alpha)$ such that $$w_{i+1}\in Son(w_{i})$$ for $i=1,\ldots,t-1.$  We will write the chain $C$ in the following form $$w_{1}- \ldots - w_{t}.$$
%\end{definition}

\subsection{ Special Chains}

The fact that a hyperlattice isomorphism preserves the parenthood, i.e., the image of a son of a given hypertuple is the son of the image of the hypertuple, leads naturally to introduce the following definition.

\begin{definition}
 A chain $C:w_{1}-w_{2}-\dots-w_{t}$ in $V(\alpha)$ is called {\em special} if it is of  maximal length satisfying the property that $Son(w_{i})=\{w_{i+1}\}, \ i=1, \ldots, t-1$. 
\end{definition}

We focus on special chains whose last element is zero. This subsection is devoted to determine the different types of special chains a hyperlattice may have. The results obtained are essentially consequence of Theorem \ref{teo:son1}.  

Observe that $V(\alpha_{1})$ is  a special chain ending at zero. 
We analyze next $V(\alpha_{1},\alpha_{2})$.

\begin{lemma}\label{lemma:specialr=2} 
Let $V(\alpha_{1}, \alpha_{2})$ a hyperlattice.

\begin{itemize}
\item  If $\alpha_{1}-\alpha_{2}=1$, then  $V(\alpha_{1}, \alpha_{2})$
is the special chain ending at zero: $$(\alpha_{1},\alpha_{2})-(\alpha_{1}-1,\alpha_{2})-\dots-(1,0)-(0,0).$$

\item If
$\alpha_{1}-\alpha_{2}>1$, then  there exists only two special chains ending at zero:
\begin{equation} \label{C1-r=2}
(1,1)-(1,0)-(0,0),
\end{equation}
and
\begin{equation} \label{C2-r=2}(\alpha_{1}-\alpha_{2},0)-\dots-(0,0).
\end{equation}
\end{itemize}
\end{lemma}

\begin{proof}  
If $\alpha_1-\alpha_2=1$, the conclusion is immediate.
Assume that $\alpha_1-\alpha_2>1$. Let $u=(u_{1},u_{2})\in V(\alpha_{1},\alpha_{2})$. If   $u_{1}>u_{2}>0$, then, by Theorem \ref{teo:son1},     $\Son(u)=\{(u_{1}-1,u_{2})\}$ 
if and only if $\alpha_{1}-u_{1}=\alpha_{2}-u_{2}$, and since $\alpha_{1}-\alpha_{2}>1$, 
$$\Son(\Son(u))=\{(u_{1}-2,u_{2}),(u_{1}-1,u_{2}-1)\},$$
therefore $u$ cannot be in a special chain ending at zero.

By Theorem \ref{teo:son1}, and following a similar argument, if   $u_{1}=u_{2}>0$,
 then  $(u_{1},u_{1})$ is in a special chain ending at zero if and only if $u_{1}=1$, and the special chain is (\ref{C1-r=2}).

Finally,  if $u_2=0$, then $(u_{1},0)$ is always in the special chain ending at zero (\ref{C2-r=2}).
%$(\alpha_{1}-\alpha_{2},0)-( \alpha_{1}-})\alpha_{2}-1,0)- \dots -(0,0)$.
\end{proof}

In the following results $k$ and $q$ are defined as in Theorem \ref{teo:son1}.

\begin{lemma}\label{lemm:notspecial} Let $V(\alpha_{1},\ldots,\alpha_{r})$ be a hyperlattice  with $r\geq 3$.  A hypertuple  of the form $u=(u_{k},\ldots,u_{k},u_{k+1},\ldots,u_{q},0,\ldots,0)$, $1<k<q$, cannot be in any special chain ending at zero.
\end{lemma}
\begin{proof} By Theorem~\ref{teo:son1}, $\card(\Son(u))=1$ if and only if $\alpha_{k}-u_{k}=\alpha_{k+1}-u_{k+1}=\dots =\alpha_q-u_q$, and in that  case $$\Son(u)=\{(u_1,\ldots, u_k-1, u_{k+1}, \ldots, u_q, 0 \ldots,0)\}.$$
Observe that  $\card(\Son(\Son(u))=1$ if and only if $\alpha_{k-1}-u_{k-1}=\alpha_k-(u_k-1)=\alpha_{k+1}-u_{k+1}=\dots =\alpha_q-u_q$, which is not true. Therefore, the unique-son chain stops at $\Son(u)$ and does not end at zero.
\end{proof}

\begin{lemma}\label{q>2}
Let $V(\alpha_{1},\ldots,\alpha_{r})$, $r\geq 2$. Then,  
\begin{equation}\label{C1-qk}(1,1,\ldots,1)-(1,\ldots,1,0)-\dots-(0,\ldots,0).
\end{equation}
is a special chain ending at zero.
\end{lemma}

\begin{proof}
It is immediate that the chain (\ref{C1-qk}) is a unique-son chain. It cannot be enlarged because the hypertuple $(1,1,\ldots,1)$ has a brother
$$\Son(2, 1, \dots, 1)=\{(1,1,\ldots,1), (2, 1, \dots, 1,0)\}.$$

\end{proof}

%\begin{lemma}\label{q>2}
%Let $V(\alpha_{1},\ldots,\alpha_{r})$, $r\geq 3$, be a hypertlattice. A hypertuple  $u=(u_{1},\ldots,u_{q},0,\ldots,0)\in V(\alpha_{1},\ldots,\alpha_{r})$, $q\geq3$,   belongs to a special chain ending at zero only if $u_1=\dots=u_q=1$, and the special chain is 
%\begin{equation}\label{C1-qk}(1,1,\ldots,1)-(1,\ldots,1,0)-\dots-(0,\ldots,0).
%\end{equation}
%\end{lemma}

%\begin{proof}
%As in Theorem \ref{teo:son1}, let $k=\max\{i:u_1=\dots=u_i\}$.
 
%If $1<k<q$, by Lemma \ref{lemm:notspecial}  the hypertuple $u$ cannot belong to any special chain ending at zero, therefore either $k=1$ or $k=q$. 

%If $k=q$, then $\Son(u)=\{\bar u=(u_1,\ldots, u_1, u_1-1,0,  \stackrel{(r-k)}{\ldots}, 0)\}$. 

%If $u_1-1>0$, then 
%$\card(\Son(\bar u))=1$ if and only if $\alpha_{q-1}-u_{1}=\alpha_q-(u_1-1)$, and  in that  case $\Son(\bar u)=\{\bar{\bar u}=(u_1,\ldots, u_1, u_{1}-1, u_1-1, 0 \ldots,0)\}$. 
%In turn, $\card(\Son(\bar{\bar u}))=1$  if and only if $\alpha_{q-2}-u_{1}=\alpha_{q-1}-(u_{1}-1)=\alpha_q-(u_1-1)$, which cannot occur. Hence, the unique-son chain does not end at zero.

%If $u_1-1=0$, then $\Son(u)=\{\bar u=(u_1,\ldots, u_1, 0, \stackrel{(r-k+1)}{\ldots}, 0)\}$.

%Repeating the argument, the  special chain ending at zero is (\ref{C1-qk}).
%\end{proof}

\begin{lemma} \label{q=12}
Let $V(\alpha_{1},\ldots,\alpha_{r})$, $r\geq 3$, be a hyperlattice. A hypertuple  $u=(u_{1},\ldots,u_{q},0,\ldots,0)\in V(\alpha_{1},\ldots,\alpha_{r})$, $q\leq2$,   belongs to a special chain ending at zero only in the following cases: 
\begin{itemize}
\item If $q=1$, then $u$ belongs to the  special chain ending at zero 
\begin{equation}\label{C1-q1}
(\alpha_1-\alpha_2, 0 \ldots,0)-(\alpha_1-\alpha_2-1, 0 \ldots,0), \cdots -(1, 0 \ldots,0)-(0, 0 \ldots,0).
\end{equation}

\item If $q=2$ and $\alpha_1-\alpha_2>1$,  then $u=(1,1,0, \ldots,0)$ and the special chain ending at zero is (\ref{C1-qk}).

%If, additionally,  $\alpha_1-\alpha_2=1$, then $u$ also belongs to (\ref{q2}).

\item If $q=2$ and  $\alpha_1-\alpha_2=1$, then  $u$ belongs to the special chain ending at zero 
\begin{equation}\label{q2}
\begin{alignedat}{2}(\alpha_2-\alpha_3+1,\alpha_2-\alpha_3, 0, \ldots,0)-(\alpha_2-\alpha_3,\alpha_2-\alpha_3, 0,  \ldots,0)- \dots \\
\hspace{4cm}  \dots - (1,0, 0, \ldots,0)-(0,0,0,  \ldots,0),
\end{alignedat}
\end{equation} and if 
$u=(1,1,0, \ldots,0)$, then it also belongs to (\ref{C1-qk}).
\end{itemize}
\end{lemma}

\begin{proof}
 We analyze the possible special chains ending at zero containing $u=(u_{1},\ldots,u_{q},0,\ldots,0)$ depending on the value of $q$.
%\begin{itemize}
%\item 

If $q=1$ it is easy to see that 
$$(\alpha_1-\alpha_2, 0 \ldots,0)-(\alpha_1-\alpha_2-1, 0 \ldots,0), \cdots -(1, 0 \ldots,0)-(0, 0 \ldots,0),$$ is a special chain  ending at zero, and $u$ belongs to it.

%\item 

Let $q=2$. 
Clearly, $u=(1,1,0,\ldots,0)$ 
belongs to (\ref{C1-qk}).

If $u_1=u_2>1$, 
then $\Son(u)=\{\bar u=(u_2,u_2-1, 0 \ldots,0)\}$. Moreover, 
$\card(\Son(\bar u))=1$ if and only if $\alpha_1-u_1=\alpha_2-(u_2-1)$ if and only if  $\alpha_1=\alpha_2+1$, and in such a case $\Son(\bar u)=\{\bar{\bar u}=(u_2-1,u_2-1, 0 \ldots,0)\}$. Repeating the argument, we get
$$(u_2,u_2, 0 \ldots,0)-(u_2,u_2-1, 0 \ldots,0)-(u_2-1,u_2-1, 0 \ldots,0)- \cdots $$
$$ \cdots -(1,0, 0 \ldots,0)-(0,0, 0 \ldots,0).$$
As $\alpha_2-u_2\geq \alpha_3$, the longest chain is obtained taking $u_2=\alpha_2-\alpha_3$, and as $\alpha_1=\alpha_2+1$, we obtain the chain (\ref{q2}).

If  $u_1>u_2$, then $\Son(u)=\{\bar u=(u_1-1,u_2, 0 \ldots,0)\}$ if and only if $\alpha_1-u_1=\alpha_2-u_2$. In turn,  $\card(\Son(\bar u))=1$ if and only if  $\alpha_1-(u_1-1)=\alpha_2-u_2$, which is not true. Hence,  $u_1-1=u_2$ and $\Son(\bar u)=\{(u_1-1,u_2-1, 0 \ldots,0)\}$. Notice that it implies that $\alpha_1=\alpha_2+1$, and this case is the same as the previous one.
%and it is of length $2(\alpha_1-\alpha_3)$. 
%\end{itemize}
\end{proof}

The analysis performed on special chains is summarized in the next theorem.

\begin{theorem}\label{theo:special}
In the hyperlattice  $V(\alpha_{1},\ldots,\alpha_{r})$, $r\geq 2$, there exist at most two of the three possible types of special  chains ending at zero,  $C_1, C_2$  and $C_3$,  of lengths $r+1, \alpha_1-\alpha_2+1$ and $2(\alpha_1-\alpha_3)$, respectively,  defined as 
follows:

$$\hspace{-1cm}\begin{array}{cccc}
& C_1 & C_2 \ (\text{if } \alpha_{1}-\alpha_{2}>1)  &  C_3 \ (\text{if } \alpha_{1}-\alpha_{2}=1) \\
\hline
& \\
& (1, 1, \ldots, 1, 1) & (\alpha_1-\alpha_2,  0, \ldots, 0)  & (\alpha_1-\alpha_3, \alpha_2-\alpha_3, 0, \ldots, 0) \\
&|&|&|\\

&  (1, 1, \ldots, 1, 0) & (\alpha_1-\alpha_2-1, 0, \ldots, 0)  & (\alpha_1-\alpha_3-1, \    \alpha_2-\alpha_3, 0,\ldots, 0) \\
&|&|&|\\
 & & & (\alpha_1-\alpha_3-1, \alpha_2-\alpha_3-1, 0, \ldots, 0) \\
&   \vdots & \vdots & \vdots \\
&|&|&|\\
&   (1, 0, \ldots, 0, 0) & (1, 0, \ldots, 0) &  (1, 0, \ldots, 0)\\
&|&|&|\\

&   (0, 0, \ldots, 0, 0) & (0, 0,  \ldots, 0) & (0, 0, \ldots, 0) \\
& \\
%{\color{blue} length\, : &  r+1, & \alpha_1-\alpha_2+1, & 2(\alpha_{1}-\alpha_3). }\\
\end{array} $$

\end{theorem}

\begin{proof}
Observe that, as already seen,  $V(\alpha_1)$ is trivially a special chain ending at zero. For $r\geq 2$, as a consequence of  Lemmas \ref{lemma:specialr=2} to \ref{q=12}, if $\alpha_1-\alpha_2>1$ the only special chains in $V(\alpha)$ are of type $C_1$ and $C_2$, and  if $\alpha_1-\alpha_2=1$, they are of type $C_1$ and $C_3$ (we take $\alpha_3=0$ when $r=2$).
\end{proof}

\begin{remark}\label{rem:existenciaespeciales} As mentioned, in a lattice there exist at most  two special chains ending at $0$. More precisely, for $r=1$ there exists a single special chain $C_2$ (observe that in this case the chain of type $C_1$ is reduced to $1-0$ and is included in $C_2$, therefore, $1-0$ is not a special chain). If $r=2$ and $\alpha_1-\alpha_2=1$, then there exists only the special chain $C_3$ (now $C_3$ includes the chains of type $C_2$ and $C_1$).
In any other case, there exist two special chains: $C_1$ and $C_2$ if $r\geq 2$ and $\alpha_1-\alpha_2>1$, or $C_1$ and $C_3$  if $r\geq 3$ and $\alpha_1-\alpha_2=1$ (in the last case  $C_1$ and  $C_3$ contain the chain of type $C_2$).
\end{remark}

\begin{example} 

\begin{itemize} \
\item 
Example of a chain of type $C_2$: Let $ V(7,3,2,1)$. Then,
%$\alpha=(7,3,2,1)$, ,  
$C: (4,0,0,0)-(3,0,0,0) -(2,0,0,0)-(1,0,0,0)-(0,0,0,0)$ is a special chain of type $C_2$, 
% in $V(7,3,2,1)$, 
of $\length(C)=\alpha_1-\alpha_2+1=5$.

\item 
Example of a chain of type $C_3$: 
%, $k=2$, $\alpha_1=\alpha_2+1$: {\color{blue}
%$\alpha=(5,4,3,2)$,
%$\alpha=(5,4,0,0)$, $u=(3,3,0,0)$.
%$w= (3,3,0,0)-(3,2,0,0)-(2,2,0,0)-(2,1,0,0)-(1,1,0,0)-(1,0,0,0)-(0,0,0,0)$, $\length(w)=7=2\cdot 3+1$.
%
%$q=2$, $k=1$, $\alpha_1=\alpha_2+1$: $\alpha=(5,4,2,2)$, $u=(4,3,0,0)$.
%$w= (4,3,0,0)-(3,3,0,0)-(3,2,0,0)-(2,2,0,0)-(2,1,0,0)-(1,1,0,0)-(1,0,0,0)-(0,0,0,0)$, $\length(w)=8=2\cdot 3+2$.
%
%{\color{blue} The maximal chain is }
%In both cases the special chain is
 Let $ V(6,5,3,2)$. Then,
$C: (3,2,0,0)-(2,2,0,0)-(2,1,0,0)-(1,1,0,0)-(1,0,0,0)-(0,0,0,0)$,  is a special chain of type $C_3$,
% in $V(\alpha)$, 
of $\length(C)=2 (\alpha_1-\alpha_{3})=6$.

\item Example of a chain of type $C_1$: Let $ V(5,4,3,2)$. Then,
$C: (1,1,1,1)-(1,1,1,0)-(1,1,0,0)-(1,0,0,0)-(0,0,0,0)$  is a special chain of type $C_1$, of $\length(C)=r+1=5$.
\end{itemize}

\end{example}

\subsection{Riding Special Chains}

Special chains are key tools to study isomorphisms of  hyperlattices,
but we also need another type of chains  associated to special chains.

\begin{definition}
  Given $p\in \{1,2,3\}$, a chain $RC_{p}:w_{1}-w_{2}-\dots -w_{l}$  rides on a special chain $C_{p}$ if it is of maximal length satisfying that $\exists j \in \{1, \ldots, l-1\}$ such that
  $\Son(w_{i})=\{w_{i+1}\}$ for  $1\leq i\leq j $, for $j+1\leq i\leq l-1$, $\Son(w_{i})=\{w_{i+1},w'_{i+1}\}$  with $w'_{i+1}\in C_{p}$, and $\Son(w_{l})=\{w'_{l+1}\}$ with $w'_{l+1}\in C_{p}$. We will refer to $RC_p$ as a {\em riding special chain} or as a {\em riding chain on $C_p$}.

 \end{definition}
To describe riding special chains in $V(\alpha_{1},\ldots,\alpha_{r})$ we analyze separately the cases when $r=2$, $r=3$ and $r>3$.

\begin{proposition}\label{propo:ridingCASE2b}   Let $V( \alpha_{1},\alpha_{2})$ be a hyperlattice.
\begin{enumerate}
\item Let  $\alpha_{1}-\alpha_{2}>1$. Then, 
the riding chains are:
	\begin{enumerate}
    	\item  $RC_1$:
		\begin{itemize}
		\item If $\alpha_{2}>1$, then  $RC_{1}: (2,2)-(2,1)-(2,0)$.
 
         	\item If $\alpha_{2}=1$ and $\alpha_{1}\neq 3$, then $ RC_{1}: (2,1)-(2,0)$.
      \item  If $\alpha_{2}=1$ and $\alpha_{1}=3$, then $ RC_{1}:(3,1)- (2,1)-(2,0)$.    	\end{itemize}
	\item $RC_2$:
    $(\alpha_{1}-\alpha_{2}+1,1)-(\alpha_{1}-\alpha_{2},1)-\ldots-(1,1)$.       	\end{enumerate}
\item If $\alpha_{1}-\alpha_{2}=1$, then there exists a $C_3$ chain, 
%(with $C_1\subset C_3$), 
but it does not exists any $RC_3$.

\end{enumerate}
     \end{proposition}
 \begin{proof}

 \begin{enumerate}
\item Let $\alpha_{1}-\alpha_{2}>1$.
	\begin{enumerate}
    	\item Existence of $RC_1$ chains:
		\begin{itemize}
    \item  If $\alpha_{2}>1$, then $\Son(2,2)=\{(2,1)\}$, $\Son(2,1)=\{(2,0),(1,1)\}$ and $\Son((2,0))=\{(1,0)\}$.  Moreover, observe that $(2,2)$ has a unique  father $(3,2)$ and $\Son(3,2)=\{(3,1),(2,2)\}$, therefore there exist a riding chain on $C_1$ of length 3:  
		$$\begin{array}{ccccccc}
   (3,2)\rightarrow &{\bf (2,2)} \rightarrow & {\bf (2,1)} \rightarrow & {\bf (2,0)} &\\
 \downarrow & & \downarrow & \downarrow&  \\
 (3,1) & & (1,1) \rightarrow & (1, 0) \rightarrow& (0,0) 
 \end{array}$$
\item The proof of the case $\alpha_{2}=1$ and $\alpha_{1}> 3$ is  analogous. Notice that $ \Son(3,1)=\{(2,1),(3,0)\}$:  
$$\begin{array}{ccccc}
  (3,1)\rightarrow & {\bf (2,1)} \rightarrow & {\bf (2,0)} &\\
 \downarrow &  \downarrow & \downarrow &   \\
 (3,0) & (1,1) \rightarrow &  (1, 0) \rightarrow& (0,0) 
 \end{array}$$

 \item If $\alpha_{2}=1$ and $\alpha_{1}=3$, the chain indicated cannot be enlarged, for
	$ V(3,1)=C_1\cup RC_1$:
	$$\begin{array}{ccccc}
  {\bf (3,1)}\rightarrow & {\bf (2,1)} \rightarrow & {\bf (2,0)} &\\
 & \downarrow & \downarrow &   \\
 & (1,1) \rightarrow  &  (1, 0) \rightarrow& (0,0) 
 \end{array}$$
        \end{itemize}

    \item  Existence of a $RC_2$ chain:   $\Son(\alpha_{1}-\alpha_{2}+1,1)=\{(\alpha_{1}-\alpha_{2},1)\}$, for $k=0,\ldots,\alpha_{1}-\alpha_{2}-2$ we have $\Son(\alpha_{1}-\alpha_{2}-k,1)=\{(\alpha_{1}-\alpha_{2}-(k+1),1),(\alpha_{1}-\alpha_{2}-k,0)\}$ with
 $(\alpha_{1}-\alpha_{2}-k,0)\in C_{2}$,  and for $k=\alpha_{1}-\alpha_{2}-1$, $\Son(1,1)=\{(1,0)\}$.
Moreover,  
 \begin{itemize}
 \item if $\alpha_2>1$, then $(\alpha_{1}-\alpha_{2}+1,2)\in V(\alpha)$, but $\Son(\alpha_{1}-\alpha_{2}+1,2)=\{(\alpha_{1}-\alpha_{2}+1,1), (\alpha_{1}-\alpha_{2},2)\}$:
 {\footnotesize
 $$\hspace{-1cm}\begin{array}{ccccccc}
  (\alpha_{1}-\alpha_{2}+1,2)\rightarrow &{\bf (\alpha_{1} -\alpha_{2}+1,1)} \rightarrow & {\bf (\alpha_{1}-\alpha_{2},1)}\rightarrow& \ldots &  \rightarrow {\bf (1,1)} &\\
 \downarrow & & \downarrow& \dots & \quad \downarrow & \\
 (\alpha_{1}-\alpha_{2},2) & & (\alpha_{1} -\alpha_{2}, 0)\rightarrow& \ldots & \rightarrow (1,0) & \rightarrow(0,0) 
 \end{array}$$
}
 
 \item  if $\alpha_2=1$, then $V(\alpha_1, 1)=C_2\cup RC_2$.
 \end{itemize}
    \end{enumerate}
    \item Let $\alpha_{1}-\alpha_{2}=1$. The consequence   is immediate, for $V(\alpha_1,\alpha_2)=C_3$.
 \end{enumerate}

 \end{proof}

\begin{proposition}\label{propo:ridingCASE3b} Let $V( \alpha_{1},\alpha_{2},\alpha_{3})$ be a hyperlattice. 
\begin{enumerate}
\item Let $\alpha_{1}-\alpha_{2}>1$. Then, 
the riding chains are: 
    \begin{itemize}
    \item   $RC_1$:
    	\begin{itemize}
        \item If $\alpha_{2}-\alpha_{3}>1$, then $RC_{1}: (2,1,1)-(2,1,0)-(2,0,0)$.
         \item 
        If $\alpha_{2}-\alpha_{3}=1$ and $\alpha_{3}\geq 2$, then $RC_{1}$:
             
              $\qquad (2,2,2)-(2,2,1)-(2,1,1)-(2,1,0)-(2,0,0)$. 
      \item If $\alpha_{2}-\alpha_{3}=1$ and $\alpha_{3}=1$, then  $RC_{1}$:
      
      $\qquad (2,2,1)-(2,1,1)-(2,1,0)-(2,0,0)$.
         \end{itemize}
    \item $RC_{2}: (\alpha_{1}-\alpha_{2}+1,1,0)-(\alpha_{1}-\alpha_{2},1,0)-\ldots-(1,1,0)$ is a $RC_2$ of $\length(RC_{2})=\alpha_{1}-\alpha_{2}+1$, and if $\alpha_1-\alpha_2=2$ and $\alpha_2-\alpha_3>1$, then $(\alpha_{1}-\alpha_{2},2,0)-(\alpha_{1}-\alpha_{2},1,0)-\ldots-(1,1,0)$ is also a $RC_2$ of the same length.
    \end{itemize}
\item If $\alpha_{1}-\alpha_{2}=1$, the riding chains are:
     	\begin{itemize}
    	\item $RC_{1}$:
       	\begin{itemize}
        \item If $\alpha_{2}-\alpha_{3}>1$, then $RC_{1}:  (2,1,1)-(2,1,0)$.
        \item  If $\alpha_{2}-\alpha_{3}=1$ and $\alpha_{3}\geq2$,  then $RC_{1}$: 

    $\qquad 
    (2,2,2)-(2,2,1)-(2,1,1)-(2,1,0)$. 
    
    \item If $\alpha_{2}-\alpha_{3}=1$ and $\alpha_{3}=1$, then $RC_{1}$: 
    
    $\qquad (3,2,1)-(2,2,1)-(2,1,1)-(2,1,0)$. 
    \end{itemize}
    
    \item $RC_{3}: (\alpha_{1}-\alpha_{3}+1,\alpha_{2}-\alpha_{3}+1,1)-(\alpha_{1}-\alpha_{3},\alpha_{2}-\alpha_{3}+1,1)-(\alpha_{1}-\alpha_{3},\alpha_{2}-\alpha_{3},1)-\dots-(1,1,1)$ and $\length(RC_{3})=2(\alpha_{1}-\alpha_{3})$.
   % $\alpha_{1}+\alpha_{2}-2\alpha_{3}+1$.
     \end{itemize}
 \end{enumerate}
 \end{proposition}
 \begin{proof} It is analogous to the proof of Proposition~\ref{propo:ridingCASE2b}.

 \begin{enumerate}
\item If $\alpha_{1}-\alpha_{2}>1$, 
    \begin{itemize}
    \item Existence of $RC_1$ chains: 
    	\begin{itemize}
        \item If $\alpha_{2}-\alpha_{3}>1$, then  $(2,1,1)-(2,1,0)-(2,0,0)$ is a riding chain on  $C_1$. Observe that it cannot be enlarged because  $(2,1,1)\in \Son \{(2,2,1), (3,1,1)\}$, but in both cases there are other sons ($\Son(2,2,1)=\{(2,1,1), (2,2,0)\}$, $\Son(3,1,1)=\{(2,1,1), (3,1,0)\}$), hence $\length(RC_{1})=3$:

        {\small        
        $$\hspace{-0.25cm}\begin{array}{ccccccc}
 (2,2,0) \leftarrow &  (2,2,1) & \\
 & & \searrow &\\
  & & & {\bf (2,1,1)}\rightarrow & {\bf (2,1,0) }\rightarrow & {\bf (2,0,0)} &\\
  & &\nearrow & \downarrow & \downarrow & \downarrow &  \\
 (3,1,1)\leftarrow & (3,1,0) & &(1,1,1) \rightarrow &  (1,1,0) \rightarrow & (1,0, 0) \rightarrow& (0,0,0) 
 \end{array}$$
} 
         \item If $\alpha_{2}-\alpha_{3}=1$ and  $\alpha_{3}\geq 2$, then
         it is easy to check that   the chain 
         $(2,2,2)-(2,2,1)-(2,1,1)-(2,1,0)-(2,0,0)$ is a riding  chain on   $C_1$.
         \item If $\alpha_{2}-\alpha_{3}=1$ and  $\alpha_{3}=1$, the chain $(2,2,1)-(2,1,1)-(2,1,0)-(2,0,0)$ is a riding chain on $C_1$.
         \end{itemize}
    \item  Existence of  $RC_2$ chains: the chain $(\alpha_{1}-\alpha_{2}+1,1,0)-(\alpha_{1}-\alpha_{2},1,0)-\ldots-(1,1,0)$ 
      rides on $C_2$. 
      Notice  that  it cannot be enlarged:

      If $\alpha_2-\alpha_3>1$, 
       the tuple 
$(\alpha_{1}-\alpha_{2}+1,1,0)$ has two fathers, each one of them with more than one son.
Notice also that if, additionally,   $\alpha_1-\alpha_2>2$, then $\Son{(\alpha_{1}-\alpha_{2},2,0)}=\{(\alpha_{1}-\alpha_{2},1,0), (\alpha_{1}-\alpha_{2}-1,2,0)\}$. The behavior is sketched in the next scheme:

If $\alpha_1-\alpha_2>2$ and $\alpha_2-\alpha_3>1$, 

{\footnotesize
$$\hspace{-1cm}\begin{array}{ccccccccc}
 (\alpha_{1}-\alpha_{2}+1,2,0) & \rightarrow&(\alpha_{1}-\alpha_{2},2,0)& \rightarrow & (\alpha_{1}-\alpha_{2}-1,2,0) \\
&\searrow & & \searrow \\
  & & {\bf (\alpha_{1}-\alpha_{2}+1,1,0)}& \rightarrow & {\bf (\alpha_{1}-\alpha_{2},1,0)}\rightarrow & \ldots &  \rightarrow {\bf (1,1,0)} &\\
& \nearrow &  & \nearrow & \downarrow&  & \downarrow & \\
 (\alpha_{1}-\alpha_{2}+1,1,1) &\rightarrow & (\alpha_{1}-\alpha_{2},1,1) 
& & (\alpha_{1}-\alpha_{2},0,0) \rightarrow &  \ldots & \rightarrow (1,0,0) & \rightarrow(0,0,0) \\& & &\searrow \\
& &  & & (\alpha_{1}-\alpha_{2}-1,1,1) 
 \end{array}$$
}
\medskip

 But,  if $\alpha_1-\alpha_2=2$ (and $\alpha_2-\alpha_3>1$), then we have $\Son{(\alpha_{1}-\alpha_{2},2,0)}=\{(\alpha_{1}-\alpha_{2}-1,2,0)\}$, therefore, $(\alpha_{1}-\alpha_{2},2,0)-(\alpha_{1}-\alpha_{2},1,0)-\ldots-(1,1,0)$ is also a riding chain  on $C_2$ of the same length.
 
%Ver 531, 631, 421

 If $\alpha_1-\alpha_2=2$, then  $(\alpha_1-\alpha_2+1,1,0)$ has only one father $(\alpha_1-\alpha_2+1,1,1)$, but it has two sons, therefore the riding chain remains the same. 
 
    	\end{itemize}
     
\item If $\alpha_{1}-\alpha_{2}=1$,
    	\begin{itemize}
    	\item Existence of $RC_1$ chains: the proof is immediate.
	       	
    \item Existence of a $RC_3$ chain:  albeit the proof is easy,  we include a sketch of it 
  when $\alpha_3>1$,   
{\footnotesize
$$\hspace{-1cm}\begin{array}{cccccc}
& \\
& &  (\alpha_1-\alpha_3+1, \alpha_2-\alpha_3+1, 2) \\
&\swarrow &  \downarrow &  \\
 (\alpha_1-\alpha_3, \alpha_2-\alpha_3+1, 2)&  & {\bf (\alpha_1-\alpha_3+1, \alpha_2-\alpha_3+1, 1)} & & \\
& & \downarrow & \\
& & {\bf (\alpha_1-\alpha_3, \alpha_2-\alpha_3+1, 1)} & \\
& & \downarrow & \\
& & {\bf (\alpha_1-\alpha_3, \alpha_2-\alpha_3, 1}) & \rightarrow & (\alpha_1-\alpha_3, \alpha_2-\alpha_3, 0) \\
& & \downarrow & &\downarrow \\
& & \vdots & &  \vdots \\
& & \downarrow & & \downarrow \\
& &  {\bf (1, 1, 1)} & \rightarrow & (1,1, 0) \\
& & & & \downarrow \\
& & & &  (1,0, 0) \\
& & & & \downarrow \\
& & & & (0, 0, 0) \\

%{\color{blue} length\, : &  r+1, & \alpha_1-\alpha_2+1, & 2(\alpha_{1}-\alpha_3). }\\
\end{array} $$
}
\noindent
If $\alpha_3=1$, then $V(\alpha_1, \alpha_2, \alpha_3)=C_3\cup RC_3$.
   \end{itemize}
 \end{enumerate}

 \end{proof}

 The riding chains when $r>3$ are described in the following proposition.
\begin{proposition}\label{propo:ridinggeneral} Let $V(\alpha_{1},\ldots,\alpha_{r})$ be a hyperlattice with $r>3$. 

\begin{enumerate}
\item If $\alpha_{1}-\alpha_{2}>1$, there exists   a $RC_1$ of $\length(RC_1)=r$:
\begin{equation} \label{rc1_r4_1}
 (2,1,\ldots,1,1)-(2,1,\ldots,1,0)-\ldots-(2,0,\ldots,0),
\end{equation}
and  a $RC_2$ of $\length(RC_2)=\alpha_{1}-\alpha_{2}+1$:
\begin{equation} \label{rc2_r4_1}
 (\alpha_{1}-\alpha_{2}+1,1,0,\ldots,0)-(\alpha_{1}-\alpha_{2},1,0,\ldots,0)-\ldots-(1,1,0,\ldots,0).
\end{equation}

\item  If $\alpha_{1}-\alpha_{2}=1$, there exists  an $RC_1$ of $\length(RC_1)=r-1$:

\begin{equation}\label{rc1_r4_2}
 (2,1,\ldots,1,1)-(2,1,\ldots,1,0)-\ldots-(2,1,0,\ldots,0),
\end{equation}
and  a $RC_3$ of $\length(RC_3)= 2(\alpha_{1}-\alpha_{3})$:
\begin{equation} \label{rc3_r4_2}
\begin{alignedat}{2}
 (\alpha_{1}-\alpha_{3}+1,\alpha_{2}-\alpha_{3}+1,1,0, 0,\ldots,0)- \hspace{4cm}\\(\alpha_{1}-\alpha_{3},\alpha_{2}-\alpha_{3}+1,1,0,\ldots,0)- %\ldots 
 \ldots -(1,1,1,0,\ldots,0).
\end{alignedat}
\end{equation}
\end{enumerate}
\end{proposition}

\begin{proof}

\begin{enumerate}
    \item 
 Let $\alpha_{1}-\alpha_{2}>1$. Observe that for $i=0, \ldots, r-2$, 
\begin{equation} \label{rc1_r4_sons}
\begin{alignedat}{2}
\Son(2,1, \stackrel{(r-i-1)}{\ldots},1,0,\stackrel{(i)}{\ldots},0)=\hspace{2cm}\\
\{(2,1,\stackrel{(r-i-2)}{\ldots},1,0,\stackrel{(i+1)}{\ldots},0),(1,1,\stackrel{(r-i-1)}{\ldots},1,0,\stackrel{(i)}{\ldots},0)\},
\end{alignedat}
\end{equation}
and $$\Son(2,0,\ldots,0)=\{(1,0,\ldots,0)\}.$$ 
Moreover, 
$$(2,1,\ldots,1)\in \Son(2,2,1,\ldots,1)\cap \Son(3,1,\ldots,1),$$ but in both cases it is not the unique son. Therefore, the chain (\ref{rc1_r4_1}) is a $RC_{1}$ in $V(\alpha_{1},\ldots,\alpha_{r})$
 of length $r$. A sketch of the riding chain appears in the next scheme:

%{\color{red} OJO! $(\alpha_{1}-\alpha_{2}+1,2,0,\ldots,0)$ puede no pertenecer a la lattice: 
%$\alpha=(5,3,2,1), RC_2: (3,1,0,0)-(2,1,0,0)-(1,1,0,0)$, $(3,1,0,0)\in \Son(3,1,1,0)$, pero $(3,2,0,0)\notin V(5,3,2,1)$}

{\footnotesize
$$\hspace{-2cm}\begin{array}{ccccccccc}
& \\
%(3,1,\ldots,1,0) & & &&(2,2, 1, \ldots, 1,0)\\
% \uparrow & & & &\uparrow &\\
& & (3,1,\ldots,1) & &  & & (2,2, 1, \ldots, 1)  & \\
& \swarrow && \searrow & & \swarrow && \searrow\\
(3,1,\ldots,1,0) & & & & {\bf (2,1,\ldots,1)} & \rightarrow & (1, 1, \ldots, 1)  & & (2,2, 1, \ldots, 1,0)\\
& & & &\downarrow & &\downarrow \\
& & & & {\bf (2, 1, \ldots, 1, 0)} & \rightarrow &(1, 1,\ldots, 1,0) \\
& & & & \downarrow & &\downarrow \\
& & & &  \vdots & \vdots & \vdots \\
& & & & \downarrow & &\downarrow \\
& & & & {\bf (2,1,0,\ldots, 0)} & \rightarrow & (1,1,0, \ldots, 0) & \\
& & & & \downarrow & &\downarrow \\
& & & & {\bf (2, 0, \ldots, 0)} & \rightarrow & (1, 0, \ldots, 0) & \\
& & & & & & \downarrow \\
& & & & & & (0, 0, \ldots,  0) &  \\

%{\color{blue} length\, : &  r+1, & \alpha_1-\alpha_2+1, & 2(\alpha_{1}-\alpha_3). }\\
\end{array} $$
}
Concerning the chain (\ref{rc2_r4_1}), observe that:
$$\Son(\alpha_{1}-\alpha_{2}+1,1,0,\ldots,0)=\{(\alpha_{1}-\alpha_{2},1,0\ldots,0)\},$$
and for $k=0,\ldots,\alpha_{1}-\alpha_{2}-1$, $$\Son(\alpha_{1}-\alpha_{2}-k,1,0,\ldots,0)=\{(\alpha_{1}-\alpha_{2}-k-1,1,0,\ldots,0),(\alpha_{1}-\alpha_{2}-k,0,0,\ldots,0)\}.$$
Moreover, if $\alpha_2-\alpha_3>1$, then
$$(\alpha_{1}-\alpha_{2}+1,1,0,\ldots,0)\in \Son(\alpha_{1}-\alpha_{2}+1,2,0,\ldots,0)\cap \Son(\alpha_{1}-\alpha_{2}+1,1,1,0,\ldots,0).$$
Then, the chain (\ref{rc2_r4_1}) is a  $RC_{2}$ in $V(\alpha_{1},\ldots,\alpha_{r})$ of length $\alpha_{1}-\alpha_{2}+1$.

If $\alpha_2-\alpha_3=1$, then $(\alpha_{1}-\alpha_{2}+1,1,0,\ldots,0)$ has a single father $(\alpha_{1}-\alpha_{2}+1,1,1,0, \ldots,0)$, but 
$$\Son(\alpha_{1}-\alpha_{2}+1,1,1,0, \ldots,0)=\{(\alpha_{1}-\alpha_{2}+1,1,1,0, \ldots,0), (\alpha_{1}-\alpha_{2}+1,1,0,\ldots,0)\},$$
therefore, $(\alpha_{1}-\alpha_{2}+1,1,1,0, \ldots,0)$ does not belong to the riding chain. Hence, the chain (\ref{rc2_r4_1}) is a $RC_2$ in $V(\alpha_{1},\ldots,\alpha_{r})$ of length $\alpha_{1}-\alpha_{2}+1$. A scheme of the behaviour appears next:

%{\color{blue}
%Concerning (\ref{rc2-r3_1}), $(\alpha_{1}-\alpha_{2}+1,1,0,\ldots,0)\in V(\alpha_{1},\ldots,\alpha_{r})$, and $$\Son(\alpha_{1}-\alpha_{2}+1,1,0,\ldots,0)=\{(\alpha_{1}-\alpha_{2},1,0\ldots,0)\}.$$For $k=0,\ldots,\alpha_{1}-\alpha_{2}-1$, $$\Son(\alpha_{1}-\alpha_{2}-k,1,0,\ldots,0)=\{(\alpha_{1}-\alpha_{2}-k-1,1,0,\ldots,0),(\alpha_{1}-\alpha_{2}-k,0,0,\ldots,0)\}.$$Moreover, since $$(\alpha_{1}-\alpha_{2}+1,1,0,\ldots,0)\in Son(\alpha_{1}-\alpha_{2}+1,2,0,\ldots,0)\cap Son(\alpha_{1}-\alpha_{2}+1,1,1,0,\ldots,0)$$the conditions for the chain (\ref{cr2-r3_2}) to be a $RC_{2}$ hold.\\}

{\footnotesize
$$\hspace{-1cm}\begin{array}{cccccc}
& \\
(\alpha_{1}-\alpha_{2},2,0,\ldots,0) & & & &(\alpha_{1}-\alpha_{2}, 1,1,0, \ldots,0)\\
 \uparrow & & & &\uparrow &\\
(\alpha_{1}-\alpha_{2}+1,2, 0, \ldots,0) & &  & & (\alpha_{1}-\alpha_{2}+1,1, 1,0, \ldots, 0)  & \\
& \searrow & & \swarrow &\\
& & {\bf (\alpha_{1}-\alpha_{2}+1,1,0,\ldots,0)}& \rightarrow & (\alpha_{1}-\alpha_{2}, 0, \ldots, 0)  & \\
& &\downarrow & &\downarrow \\
& & {\bf (\alpha_{1}-\alpha_{2}, 1,0, \ldots, 0)} & \rightarrow &(\alpha_{1}-\alpha_{2}-1,0, \ldots, 0) \\
& & \downarrow & &\downarrow \\
& &  \vdots & \vdots & \vdots \\
& & \downarrow & &\downarrow \\
& & {\bf (1,1,0,\ldots, 0)} & \rightarrow & (1,0, \ldots, 0) & \\

& & & & \downarrow \\
& & & & (0, 0, \ldots,  0) &  \\

%{\color{blue} length\, : &  r+1, & \alpha_1-\alpha_2+1, & 2(\alpha_{1}-\alpha_3). }\\
\end{array} $$
}
%{\color{red} OJO! $(\alpha_{1}-\alpha_{2}+1,2,0,\ldots,0)$ puede no pertenecer a la lattice: $\alpha=(5,3,2,1), RC_2: (3,1,0,0)-(2,1,0,0)-(1,1,0,0)$, $(3,1,0,0)\in \Son(3,1,1,0)$, pero $(3,2,0,0)\notin V(5,3,2,1)$}

\item 
Assume that $\alpha_{1}-\alpha_{2}=1$.  Observe that for $i=0, \ldots, r-3$, (\ref{rc1_r4_sons}) is satisfied, 
and
$$\Son(2,1, 0,\ldots,0)=\{(1,1,0,\ldots,0)\}.$$
Since 
$$(2,1,\ldots,1), (2,2,1,\ldots,1,0)\in \Son(2,2,1,\ldots,1),$$  the chain (\ref{rc1_r4_2}) is a $RC_{1}$ in $V(\alpha_{1},\ldots,\alpha_{r})$
 of length $r-1$.\\

 On the other hand, 
 $$\Son(\alpha_{1}-\alpha_{3}+1,\alpha_{2}-\alpha_{3}+1,1,0,\ldots,0)=\{(\alpha_{1}-\alpha_{3},\alpha_{2}-\alpha_{3}+1,1,0,\ldots,0)\},$$
$$\Son(\alpha_{1}-\alpha_{3},\alpha_{2}-\alpha_{3}+1,1,0,\ldots,0)=\{(\alpha_{1}-\alpha_{3},\alpha_{2}-\alpha_{3},1,0,\ldots,0)\},$$for $k=1, \ldots, \alpha_1-\alpha_3-1$,
 $$\Son(k+1,k+1,1,0,\ldots,0)=\{(k+1,k,1,0,\ldots,0), (k+1,k+1,0,\ldots,0)\},$$
 %for $k=2, \ldots, \alpha_1-\alpha_3$,
$$\Son(k+1,k,1,0,\ldots,0)=\{(k,k,1,0,\ldots,0), (k+1,k,0,\ldots,0)\},$$
$$\Son(1,1,1,0,\ldots,0)=\{(1,1,0,\ldots,0),$$
and 
$$\{(\alpha_{1}-\alpha_{3}+1,\alpha_{2}-\alpha_{3}+1,1,0,\ldots,0),(\alpha_{1}-\alpha_{3},\alpha_{2}-\alpha_{3}+1,1,1,\ldots,0)\}=  $$
$$\Son((\alpha_{1}-\alpha_{3}+1,\alpha_{2}-\alpha_{3}+1,1,1,0,\ldots,0)),$$
therefore, (\ref{rc3_r4_2})  is a $RC_3$ in $V(\alpha_1, \ldots, \alpha_r)$ of length $2(\alpha_1-\alpha_3)$.

 {\footnotesize 
$$\hspace{-1cm}\begin{array}{ccccc}
& & RC_3  &  & C_3 \\ % (r\geq 2, \alpha_{1}-\alpha_{2}=1) \\
\\
  (\alpha_{1}-\alpha_{3}, \alpha_{2}-\alpha_{3},1,0,\ldots,0) & &    (\alpha_{1}-\alpha_{3}, \alpha_{1}-\alpha_{3},0,\ldots,0)  & &   (\alpha_{2}-\alpha_{3}, \alpha_{2}-\alpha_{3},2,0,\ldots,0)\\
\uparrow  &  &\uparrow & &\uparrow  \\
 (\alpha_{1}-\alpha_{3}, \alpha_{2}-\alpha_{3},1,1,0,\ldots,0) & & (\alpha_{1}-\alpha_{3}, \alpha_{1}-\alpha_{3},1,0,\ldots,0) &  &  (\alpha_{1}-\alpha_{3}, \alpha_{2}-\alpha_{3},2,0,\ldots,0) \\
  &\searrow &  \downarrow & \swarrow &(\text{ if } \alpha_3-\alpha_4>1)\\
%\hline
& \\
&& {\bf (\alpha_{1}-\alpha_{3},\alpha_{2}-\alpha_{3},1,0, \ldots, 0)} & \rightarrow &  (\alpha_1-\alpha_3, \alpha_2-\alpha_3, \ldots, 0) \\
&&\downarrow & & \downarrow \\
& & {\bf (\alpha_{2}-\alpha_{3},\alpha_{2}-\alpha_{3},1,0,\ldots,0)} &  \rightarrow & (\alpha_2-\alpha_3, \    \alpha_2-\alpha_3, \ldots, 0) \\
&&\downarrow & & \downarrow  \\
 && {\bf (\alpha_{2}-\alpha_{3},\alpha_{2}-\alpha_{3}-1,1,0,\ldots,0)}&  \rightarrow & (\alpha_2-\alpha_3, \alpha_2-\alpha_3-1, \ldots, 0) \\
 && \downarrow & & \downarrow \\
&&   \vdots & & \vdots \\
&& \downarrow & & \downarrow \\
& & {\bf (1, 0, \ldots, 0, 0)} &\rightarrow & (1, 0, \ldots, 0)\\
&& && \downarrow  \\
&& & & (0, 0, \ldots, 0) \\
& \\
\text{ lengths: } && \alpha_1-\alpha_2+1 & & 2(\alpha_{1}-\alpha_3) \\
\end{array} $$ 
}

 \end{enumerate}
\end{proof}

\begin{remark} This section could have been developed studying the hypertuple's relation of parenthood, due to the property of selfduality of hyperinvariant lattices (\cite{FillHL77}). In that case, the components of a special chain will have only one parent, and they  will end at $u= (\alpha_{1},\ldots,\alpha_{r})$.  It is left to the reader to state the results with the parenthood relation.\end{remark}

\section{Isomorphisms of hyperinvariant lattices}
\label{isomorphisms}

Now we have  the tools  to characterize isomorphisms of hyperlattices. As the lengths of the special chains and the riding chains must be preserved by isomorphisms, we will use them  to disregard cases that do not meet this property.

We study  first isomorphisms $f:V(\alpha_{1},\ldots,\alpha_{r}) \longrightarrow V(\beta_{1},\ldots,\beta_{s})$ with $r\neq s$, and afterwards when $r=s$.
When analyzing the special chains or the riding special chains, we will denote by $C_i, RC_i$  those in $V(\alpha_{1},\ldots,\alpha_{r})$ and by $C_i', RC_i'$  those in $V(\beta_{1},\ldots,\beta_{s})$, for $i=1,2,3$.
\subsection{Isomorphic hyperinvariant lattices when $r\neq s$.} \label{segrerneqs}

This section is essentially devoted to study isomorphisms of hyperlattices when $r\neq s$.
We start with the  following auxiliary lemma.

\begin{lemma}\label{lemC1} There are no isomorphisms $f:V(\alpha_{1},\ldots,\alpha_{r}) \longrightarrow V(\beta_{1},\ldots,\beta_{s})$ with $r\neq s$ such that $f(C_{1})=C'_{1}$.
\end{lemma}
\begin{proof} $f(C_{1})=C'_{1}$ implies $r=s$, which is a contradiction.
\end{proof}

When $r=1$ we obtain the next result.
 \begin{proposition}\label{propo:r1}
 $V(\alpha_{1})\simeq V(\beta_{1},\ldots,\beta_{s}$) if and only if $s=1$ and $\alpha_{1}=\beta_{1}$ or $s=2$, $\alpha_{1}=2\beta_{1}-1$ and $\beta_{2}=\beta_{1}-1$.
 \end{proposition}
 \begin{proof}
 If $s=1$ the result is obvious. 
If $s=2$, as $V(\beta_{1},\beta_{2})$ must have only one special chain, $\beta_{1}-\beta_{2}=1$. Since $\alpha_{1}=\beta_{1}+\beta_{2}$, from Remark \ref{rem:sum-card}  the conclusion is immediate. If $s>2$, then $V(\beta_{1},\ldots,\beta_{s})$ has always two special chains, therefore $V(\alpha_{1})$ and $V(\beta_{1},\ldots,\beta_{s})$ cannot be isomorphic.

The sufficiency is immediate (see Remark \ref{rem21}). 
 \end{proof}

 The following proposition can be proved without using riding chains, but with  them the proof  becomes easier. 
\begin{proposition}\label{propo:case23}
 $V(\alpha_{1},\alpha_{2})\simeq V(\beta_{1},\ldots,\beta_{s})$, $s\geq 3$,  %$\Leftrightarrow$ 
if and only if $\alpha=(5,2)$ and $\beta=(4,2,1)$.
\end{proposition}
\begin{proof}

Let $f: V(\alpha_{1},\alpha_{2})\rightarrow V(\beta_{1},\ldots,\beta_{s})$ be an isomorphism.

If $\alpha_{1}-\alpha_{2}=1$ then, as a consequence of Proposition \ref{propo:r1}, $V(\alpha_{1},\alpha_{2})\simeq V(2\alpha_{1}-1)$, therefore we can assume that $\alpha_{1}-\alpha_{2}>1$.

 If $\beta_{1}-\beta_{2}=1$ then, by Lemma \ref{lemC1},   $f(C_{1})=C'_{3}$. It implies that $3=2(\beta_1-\beta_3)$, which cannot occur. Hence, $\beta_{1}-\beta_{2}>1$.

    \begin{itemize}
    \item 
Assume first that $s=3$.

Taking into account Lemma \ref{lemC1},  
$f(C_{2})=C'_{1}$  and  $f(C_{1})=C'_{2}$. It means that $\alpha_{1}-\alpha_{2}=3$ and $\beta_{1}-\beta_{2}=2$.
As $f(RC_{2})=RC'_{1}$, by  Propositions~\ref{propo:ridingCASE2b} and \ref{propo:ridingCASE3b}  necessarily  $4=\alpha_1-\alpha_2+1=\length(RC_2)=\length(RC_1')$, from where we conclude that $\beta_{2}-\beta_{3}=\beta_3=1$. 
Therefore, 
$\beta=(4,2,1)$. 
The dimension condition (see  Remark \ref{rem:sum-card}) leads to  
$2\alpha_{2}+3=7$, i.e.,
$\alpha=(5,2)$.

\item Assume now that $s>3$. We prove that $f$ cannot exist. Indeed, by
Lemma \ref{lemC1}, 
$f(C_{2})=C'_{1}$, hence $\alpha_{1}-\alpha_{2}=s>3$. 
%and $\beta_{1}-\beta_{2}=2$. 
Moreover, as $f(RC_{2})=RC'_{1}$,  by Propositions \ref{propo:ridingCASE2b} and \ref{propo:ridinggeneral}, %$\length(RC_1)\leq 3< \length(RC_2')$.
$\length(RC_2)=\alpha_{1}-\alpha_{2}+1=\length(RC_1')=s$, which is a contradiction.
    \end{itemize}
The sufficiency is immediate.
\end{proof}

\begin{proposition}~\label{prop:3s}
There are no isomorphisms $f: V(\alpha_{1},\ldots,\alpha_{r})\rightarrow V(\beta_{1},\ldots,\beta_{s})$ with $r\geq 3$, $s>3$ and $r\neq s$.
\end{proposition}

\begin{proof}   
By Lemma \ref{lemC1}, $f(C_{1})\neq C'_{1}$.
\begin{itemize}
\item If  $\alpha_{1}-\alpha_{2}>1$, then necessarily $f(C_{2})=C'_{1}$. It means that $\length(C_{2})=\length(C'_{1})$, i.e., $\alpha_1-\alpha_2+1=s+1$. As 
$\length(RC_{2})=\length(RC'_{1})$,  by Propositions \ref{propo:ridingCASE3b} and \ref{propo:ridinggeneral} we get $ \alpha_{1}-\alpha_{2}+1=s$ if $\beta_1-\beta_2>1$ or $ \alpha_{1}-\alpha_{2}+1=s-1$ if $\beta_1-\beta_2=1$. In both cases we have a contradiction.

\item If  $\alpha_{1}-\alpha_{2}=1$, then necessarily $f(C_{3})=C'_{1}$. The conclusion follows as in the previous case substituting $\length(C_2)$ by $\length(C_3)=2(\alpha_1-\alpha_3)$ and $\length(RC_{2})$ by $\length(RC_{3})=2(\alpha_1-\alpha_3)$.
\end{itemize}

\end{proof}

%%%%%%%%%%%%%%%%%%%%%%%%%%%%%%%%%%%%%%%%%%%%%%%%%%%%%%%%%%%%%%%%%%%%%
\subsection{Isomorphic hyperinvariant lattices when $r=s$.} \label{segrer=s}

The first case  we analyze is  $r=s=2$.  

\begin{proposition}\label{propo:case22}
$V(\alpha_{1},\alpha_{2})\simeq V(\beta_{1},\beta_{2}) \ \text{ if and only if  } \ \alpha_{i}=\beta_{i}$
\end{proposition}
\begin{proof}
 The sufficiency is immediate.

If $\alpha_{1}-\alpha_{2}=1$, then $V(\alpha_{1},\alpha_{2})\simeq V(2\alpha_{2}+1)$. Hence,  $V(\beta_{1},\beta_{2})$ shall be isomorphic to $ V(2\beta_{2}+1)$, and it implies that $\beta_{1}-\beta_{2}=1$, $\alpha_2=\beta_2$ and  $\alpha_1=\beta_1$. 
 
It remains to analyze then, the case $\alpha_{1}-\alpha_{2}>1$ and $\beta_{1}-\beta_{2}>1$:  
 \begin{itemize} 
 \item If $f(C_{2})=C'_{2}$, then $\alpha_{1}-\alpha_{2}=\beta_{1}-\beta_{2}$. As $\alpha_{1}+\alpha_{2}=\beta_{1}+\beta_{2}$, the result follows.
  \item If $f(C_{1})=C'_{2}$ and $f(C_{2})=C'_{1}$, then  $\beta_{1}-\beta_{2}=2$ and $\alpha_{1}-\alpha_{2}=2$. By the cardinality condition (Remark \ref{rem:sum-card}) we have 
 $(\alpha_{1}-\alpha_{2}+1)(\alpha_{2}+1)=(\beta_{1}-\beta_{2}+1)(\beta_{2}+1)$, 
 from where we conclude that $\alpha_{1}=\beta_{1}$ and $\alpha_{2}=\beta_{2}$.
 \end{itemize} 
\end{proof}

\begin{proposition}~\label{prop:1}
 There are no  isomorphisms $f: V(\alpha_{1},\alpha_{2},\alpha_{3})\rightarrow V(\beta_{1},\beta_{2},\beta_{3})$,  $\alpha\neq \beta$, such that $f(C_{i})= C'_{j}$ for any $i,j=1,2,3$ and $i\neq j$.

\end{proposition}
\begin{proof}

  Assume that there exists an isomorphism $f: V(\alpha_1,\alpha_{2},\alpha_{3}) \rightarrow V(\beta_{1},\beta_{2},\beta_{3})$ where $\alpha\neq \beta$ and $f(C_i)\neq C_i'$ for some $i$. We must  analyze the following cases:
    \begin{itemize}
    \item  If $\alpha_{1}-\alpha_{2}>1$ and $\beta_{1}-\beta_{2}>1$  then $f(C_{1})=C'_{2}$ and $f(C_{2})=C'_{1}$, which implies that $\alpha_{1}-\alpha_{2}=\beta_{1}-\beta_{2}=3$. Moreover,  as $f(RC_{1})=RC'_{2}$, by Proposition \ref{propo:ridingCASE3b} we conclude that  $\length(RC_1)=\length(RC_2')=\beta_1-\beta_2+1=4$ and  $\alpha_2-\alpha_3=\alpha_3=1$. Analogously, $\beta_2-\beta_3=\beta_3=1$. Hence $(\alpha_1,\alpha_2,\alpha_3)=(\beta_1,\beta_2,\beta_3)=(5, 2, 1)$, which contradicts the fact that $\alpha\neq \beta$.

    \item   If $\alpha_{1}-\alpha_{2}>1$ and $\beta_{1}-\beta_{2}=1$,  we have two possibilities:  $f(C_{1})=C'_{3}$ and $f(C_{2})=C'_{1}$ or $f(C_{1})=C'_{1}$ and $f(C_{2})=C'_ {3}$. 

 In the first case  $2(\beta_1-\beta_3)=4$ and $\alpha_{1}-\alpha_{2}=3$. As   $\length(RC_{1})=\length(RC'_{3})$, from Proposition \ref{propo:ridingCASE3b} we conclude that $\alpha_{2}-\alpha_{3}=\alpha_3=1$. Hence, $\alpha=(5,2,1)$ and $\beta=(\beta_3+2, \beta_3+1, \beta_3)$,  and from the dimension condition (Remark \ref{rem:sum-card}) we get  $5=3\beta_3$, which is not possible.

In the second option, from $f(C_{2})=C'_ {3}$ we obtain that $\alpha_{1}-\alpha_{2}+1=2(\beta_{1}-\beta_{3})$. As $\length(RC_1)=\length(RC_1')$, from Proposition \ref{propo:ridingCASE3b} we conclude that $\alpha_{2}-\alpha_{3}=\alpha_3=1$ and $\beta_{2}-\beta_{3}=1$.  Hence, $\alpha=(5,2,1)$ and $\beta=(\beta_3+2, \beta_3+1, \beta_3)$, and the conclusion follows as before. 

    \item If $\alpha_{1}-\alpha_{2}=1$ and $\beta_{1}-\beta_{2}=1$, then
$f(C_{1})=C'_{3}$, which implies that  $4=2(\beta_{1}-\beta_{3})$. As $\length(RC_1)=\length(RC_3')$, necessarily $\alpha_2-\alpha_3=1$. Hence  $\alpha=(\alpha_3+2, \alpha_3+1, \alpha_3)$ and $\beta=(\beta_3+2, \beta_3+1, \beta_3)$. By Remark \ref{rem:sum-card}  we conclude that $\alpha_i=\beta_i, i=1, 2,3$.

    \end{itemize}
 
    \end{proof}

\begin{remark}
    From the proof of the   above proposition we can see that there is  a single hyperlattice, $V(5,2,1)$, such that there exist two isomorphisms  
    $$f: V(5,2,1)\rightarrow V(5,2,1),$$ which are the identity, and another one satisfying $f(C_{1})=C_{2}$ (see Figure 1).
    
    Analogously, for the family of hyperlattices $V(\alpha_{3}+2,\alpha_{3}+1,\alpha_{3})$ for $\alpha_{3}=1,2, \ldots$, there exist two possible isomorphisms  $$f: V(\alpha_{3}+2,\alpha_{3}+1,\alpha_{3})\rightarrow V(\alpha_{3}+2,\alpha_{3}+1,\alpha_{3}),$$ the identity isomorphism,   and another one such that  $f(C_{1})=C_{3}$. 
\begin{figure}[ht]
  \centering
  \includegraphics[width=0.8\textwidth]{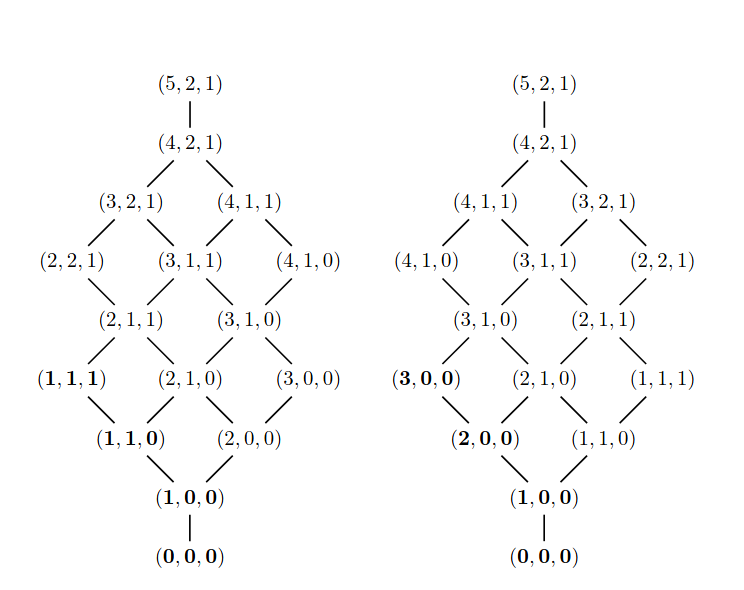}
  \caption{Hyperlattice $V(5,2,1)$}
\end{figure}

\end{remark}

\begin{proposition}\label{lemma:rCi}
 Let 
 $$f:V(\alpha_{1},\ldots,\alpha_{r})\longrightarrow V(\beta_{1},\ldots,\beta_{r})$$ be an isomorphism with $r> 3$, then  $f(C_{1})=C_{1}'$. 
 \end{proposition}
 \begin{proof}
  Assume that $f(C_{1})\neq C_{1}'$. The prove is analogous to that of Proposition \ref{prop:3s}.
 \end{proof}

 It only remains to study isomorphisms  $f: V(\alpha_{1},\ldots,\alpha_{r}) \longrightarrow V(\beta_{1},\ldots,\beta_{r})$ such that $f(C_1)=C_1'$. To do it we need some isomorphism theorems of factor lattices, which are built as quotients  by congruent relations.

%%%% Pongo aquí la congruencia ~_2

\begin{definition}
Given a hyperlattice $V(\alpha_{1},\ldots,\alpha_{r})$, we define the following congruence relations:

\begin{itemize}
\item If $\alpha_1-\alpha_2>1$ and  $r\geq 2$, we denote by $\sim_{2}$  the congruence relation defined as
$$(u_{1},\ldots,u_{r})\sim_{2}(\bar{u}_{1},\ldots,\bar{u}_{r})\Leftrightarrow u_{2}=\bar{u}_{2} ,\ldots,u_{r}=\bar{u}_{r}.$$

\item 

If $\alpha_{1}-\alpha_{2}=1$ and  $r\geq 3$,  we denote by $\sim_{3}$  the congruence relation defined as
$$(u_{1},\ldots,u_{r})\sim_{3} (\bar{u}_{1},\ldots,\bar{u}_{r})\Leftrightarrow u_{3}=\bar{u}_{3} ,\ldots,u_{r}=\bar{u}_{r}.$$

\end{itemize}

\end{definition}

The  relation $\sim_2$ allows us to build the factor lattice $V(\alpha_{1},\ldots,\alpha_{r})/\sim_{2}$.
Notice that  if $(u_{1},u_{2},u_{3}, \ldots, u_{r}) \in  V(\alpha_{1},\ldots,\alpha_{r})$, all of the elements of $ V(\alpha_{1},\ldots,\alpha_{r})$ congruent to it are 
 $$\{(u_{2},u_{2},u_{3}, \ldots, u_{r}), \ldots, (u_{2}+\alpha_{1} -\alpha_{2},u_{2},u_{3}, \ldots, u_{r})\},$$
 where $u_2\leq u_1\leq u_2+\alpha_{1} -\alpha_{2}$.  Therefore, each class contains $\alpha_1-\alpha_2+1$ elements and 
 $$\card(V(\alpha_{1},\ldots,\alpha_{r})/\sim_{2}) =(\alpha_{2}-\alpha_{3}+1)\dots(\alpha_{r}+1).$$
We  take as representative of a class in $V(\alpha_{1},\ldots,\alpha_{r})/\sim_{2}$ the first one of the previous elements; i.e., 
$$[(u_{2},u_{2},u_{3}, \ldots, u_{r})]=\{(u_{2},u_{2},u_{3}, \ldots, u_{r}), \ldots, (u_{2}+\alpha_{1} -\alpha_{2},u_{2},u_{3}, \ldots, u_{r})\}.$$
The join and meet of two classes are defined, respectively, as follows
  $$[(u_{2},u_{2},u_{3}, \ldots, u_{r})] \vee [(\bar u_{2},\bar u_{2},\bar u_{3}, \ldots, \bar u_{r})]=[(\max(u_{2},\bar u_{2}),\max(u_{2},\bar u_{2}),\ldots,\max(u_{r},\bar u_{r}))], $$
$$[(u_{2},u_{2},u_{3}, \ldots, u_{r})] \wedge [(\bar u_{2},\bar u_{2},\bar u_{3}, \ldots, \bar u_{r})]=[(\min(u_{2},\bar u_{2}),\min(u_{2},\bar u_{2}),\ldots,\min(u_{r},\bar u_{r}))]. $$

Analogously, $\sim_3$ allows us to build the factor lattice $V(\alpha_{1},\ldots,\alpha_{r})/\sim_{3}$.
Given $(u_{1},u_{2},u_{3}, \ldots, u_{r}) \in  V(\alpha_{1},\ldots,\alpha_{r})$, all of the elements of $ V(\alpha_{1},\ldots,\alpha_{r})$ congruent to it are 
 $$\{(u_{3},u_{3},u_{3},u_{4}, \ldots, u_{r}), \ldots, (u_{3}+\alpha_{1} -\alpha_{3},u_{3}+\alpha_{2}-\alpha_{3},u_{3},u_{4}, \ldots, u_{r})\},$$
 where $u_3\leq u_1\leq u_3+\alpha_{1} -\alpha_{3}$ and $u_3\leq u_2\leq u_3+\alpha_{2} -\alpha_{3}$.
Notice that each class contains $2(\alpha_1-\alpha_3)$ elements and 
$$\card(V(\alpha_{1},\ldots,\alpha_{r})/\sim_{3})= 
    (\alpha_{3}-\alpha_{4}+1)\dots(\alpha_{r}+1).$$ 
    %We will  take as representative of a class in $V(\alpha_{1},\ldots,\alpha_{r})/\sim_{3}$ the first one of the previous elements. 
    The join and meet in this quotient lattice are defined as in the previous case.

\medskip
The following lemma is satisfied.

\begin{lemma}\label{lemQuo}
We have the following isomorphisms:
\begin{enumerate}
\item
   If $r\geq 2$ and $\alpha_{1}-\alpha_{2}>1$,  then $V(\alpha_{1},\ldots,\alpha_{r})/\sim_{2} \ \simeq V(\alpha_{2},\ldots,\alpha_{r})$.
   \item If $r\geq 3$ and $\alpha_{1}-\alpha_{2}=1$,  then $V(\alpha_{1},\ldots,\alpha_{r})/\sim_{3} \ \simeq V(\alpha_{3},\ldots,\alpha_{r})$.
   \end{enumerate}
\end{lemma}
\begin{proof}
    We only prove item 1 because the proof of item 2 is analogous.
    To check that both lattices are isomorphic we prove that there is a bijective  map between the two lattices  preserving meets and joins.

    The following map  is well defined and  is clearly bijective 
    $$\begin{array}{cccc}
    \phi:&  V(\alpha_{1},\ldots,\alpha_{r})/\sim_2& \longrightarrow & V(\alpha_{2},\ldots,\alpha_{r}) \\
    & [(u_{2},u_{2},u_{3}, \ldots, u_{r})] & \mapsto & (u_{2},u_{3}, \ldots, u_{r})
    \end{array}$$

   In order to prove join preservation it is sufficient to check whether  the join of representatives of two classes maps to the join of its images, which is true. 
    The preservation of the meet can be checked  analogously.
\end{proof}

Let $V(\alpha_{1},\ldots,\alpha_{r})$ and $V(\beta_{1},\ldots,\beta_{r})$ be hyperlattices with $r\geq 2$. If $\alpha_{1}-\alpha_{2}> 1$ and $\beta_1-\beta_2>1$, we are going to prove that if there exists an isomorphism between the hyperlattices transforming $C_2$ into $C_2'$, then  $V(\alpha_{1},\ldots,\alpha_{r})/ \sim_2$ and $V(\beta_{1},\ldots,\beta_{r})/\sim_2$ are also isomorphic.

\medskip
%{\bf He encontrado que la defininicion de $\ker f$ en latices es $(u_{1},\ldots,u_{r})\in\ker f$ si y sólo si $f(u_{1},\ldots,u_{r})=(0,\ldots,0)$. Fijaos que si $f(u_{1},\ldots,u_{r})=(v_{1},\ldots,v_{r})$, entonces es automatico que $\ker \bar{f}=C_{2}$}\\

Given an isomorphism 
$$f:V(\alpha_{1},\ldots,\alpha_{r})\longrightarrow V(\beta_{1},\ldots,\beta_{r}),$$
satisfying $f(C_{2})=C'_{2}$,
and the epimorphism 
$$\pi:V(\beta_{1},\ldots,\beta_{r})\longrightarrow V(\beta_{1},\ldots,\beta_{r})/\sim_{2}, $$
let $\bar{f}=\pi \circ f$. As a consequence of Theorem~\ref{teor:iso} we have \begin{equation}\label{eqganchito2}V(\alpha_{1},\ldots,\alpha_{r})/\ker(\bar{f})\simeq V(\beta_{1},\ldots,\beta_{r})/\sim_{2}.\end{equation}

If  $\alpha_{1}-\alpha_{2}=\beta_{1}-\beta_{2}=1$, $r\geq 3$, let 
$$g:V(\alpha_{1},\ldots,\alpha_{r})\longrightarrow V(\beta_{1},\ldots,\beta_{r}),$$
be an isomorphism such that $g(C_{3})=C'_{3}$, it can be proved  in an analogous way  that 
\begin{equation}\label{eq:ganchito3}V(\alpha_{1},\ldots,\alpha_{r})/\ker(\bar{g})\simeq V(\beta_{1},\ldots,\beta_{r})/\sim_{3},
\end{equation}
where $\bar g=\pi\circ g$.

\medskip

We only need to prove the  following result.

\begin{theorem}\label{teore:iso2}
    With the above notation
    \begin{enumerate}
    \item If $r\geq 2$, $\alpha_{1}-\alpha_{2}>1$ and $\beta_{1}-\beta_{2}>1$, then
    $$V(\alpha_{1},\ldots,\alpha_{r})/\ker(\bar{f})=V(\alpha_{1},\ldots,\alpha_{r})/\sim_{2}.$$
 \item If $r\geq 3$ and $\alpha_{1}-\alpha_{2}=\beta_{1}-\beta_{2}=1$, then
    $$V(\alpha_{1},\ldots,\alpha_{r})/\ker(\bar{g})=V(\alpha_{1},\ldots,\alpha_{r})/\sim_{3}.$$
    \end{enumerate}
\end{theorem}

\begin{proof} We only prove item 1. The proof of item 2 can be done in an analogous way.
    We show that the classes of the two lattices $V(\alpha_{1},\ldots,\alpha_{r})/\ker(\bar{f})$ and $V(\alpha_{1},\ldots,\alpha_{r})/\sim_{2}$ coincide.
    
As the number of elements in a class in  $V(\beta_{1},\ldots,\beta_{r})/\sim_{2}$  is $\beta_{1}-\beta_{2}+1$, by (\ref{eqganchito2})  we know that 
$\card(V(\alpha_{1},\ldots,\alpha_{r})/\ker(\bar{f}))=\card(V(\beta_{1},\ldots,\beta_{r})/\sim_2)$, and it is  $(\beta_2-\beta_3+1)\dots (\beta_{r-1}-\beta_r+1)(\beta_r+1)$.  

By hypothesis $f(C_{2})=C'_{2}$, then $\alpha_{1}-\alpha_{2}=\beta_{1}-\beta_{2}$. As a consequence, $\card(V(\alpha_{1},\ldots,\alpha_{r})/\sim_{2})=\card(V(\alpha_{1},\ldots,\alpha_{r})/\ker(\bar{f}))$. 
Moreover, let $C$ be a class in $V(\beta_{1},\ldots,\beta_{r})/\sim_2$. Notice that the inverse image  of the class $\bar f^{-1}(C)$ contains $\alpha_1-\alpha_2+1$ elements (as does $\pi^{-1}(C)$), and if $u\in \bar f^{-1}(C)$, then $\bar f(u)=C$, i.e., the set $\bar f^{-1}(C)$ is included in the same class of $V(\alpha_{1},\ldots,\alpha_{r})/\ker(\bar{f})$. Therefore, each class in  $V(\alpha_{1},\ldots,\alpha_{r})/\ker(\bar{f})$ has $\alpha_1-\alpha_2+1$ elements. Furthermore, $C_2=\bar f^{-1}(0)$, therefore $C_2 \in V(\alpha_{1},\ldots,\alpha_{r})/\ker(\bar{f})$.

Let us see that $V(\alpha_{1},\ldots,\alpha_{r})/\sim_{2}$ and $V(\alpha_{1},\ldots,\alpha_{r})/\ker(\bar{f})$ contain the same classes.

%$$\cancel{C_{2}=[(0,0,\ldots,0)]=\{(0,\ldots,0),(1,0,\ldots,0),\ldots,(\alpha_{1}-\alpha_{2},0,\ldots,0)\}.}$$
Obviously, $C_{2} \in V(\alpha_{1},\ldots,\alpha_{r})/\sim_{2}$. From the chain $C_{2}$ we construct  classes of  $V(\alpha_{1},\ldots,\alpha_{r})/\ker(\bar{f})$, and we see that this quotient lattice contains  all of the classes of $V(\alpha_{1},\ldots,\alpha_{r})/\sim_{2}$.

Since $f(RC_{2})=RC'_{2}$, the set of elements of $RC_2$ 
   $$\{(1,1,0,\ldots,0),(2,1,0,\ldots,0),\ldots,(\alpha_{1}-\alpha_{2},1,0,\ldots,0),(\alpha_{1}-\alpha_{2}+1,1,0,\ldots,0)\}$$
    is a class in $V(\alpha_{1},\ldots,\alpha_{r})/\ker(\bar{f})$. 

Let $[(u_2,u_{2},u_{3},\ldots,u_{r})]\in V(\alpha_{1},\ldots,\alpha_{r})/\sim_2$, $u_2\neq \alpha_2$.
A father of this class is of the form $[(u_2,u_{2},\ldots, u_i+1,\ldots,u_{r})]$ for some $i\in\{3, \dots, r\}$ or
$[(u_2+1,u_{2}+1,u_{3},\ldots,u_{r})]$.

Assume that the father is $[(u_2,u_{2},\ldots, u_i+1,\ldots,u_{r})]$.  Since the element $(u_2,u_{2},\ldots, u_i+1,\ldots,u_{r})\in V(\alpha_{1},\ldots,\alpha_{r})$, 
 a  class parent  can be constructed as
 $$\begin{array}{ccc}
 (u_{2},u_{2},u_{3},\ldots,u_{i},\ldots,u_{r}) & \vee & (u_{2},u_{2},u_{3},\ldots,u_{i}+1,\ldots,u_{r})\\
 \vdots & & \vdots \\
 (u_{2}+\alpha_{1}-\alpha_{2},u_{2},u_{3}\ldots,u_{i},\ldots,u_{r})& \vee & (u_{2},u_{2},u_{3}\ldots,u_{i}+1,\ldots,u_{r})\\
 \end{array}$$
 It means that  $[(u_2,u_{2},\ldots, u_i+1,\ldots,u_{r})]\in V(\alpha_{1},\ldots,\alpha_{r})/\ker(\bar{f})$.
 
Assume now that  the father is $[(u_2+1,u_{2}+1,u_{3},\ldots,u_{r})]$. As $(u_2+1,u_{2}+1,u_{3},\ldots,u_{r})\in V(\alpha_{1},\ldots,\alpha_{r})$, an analogous construction proves that
$$\{(u_{2}+1,u_{2}+1,u_{3},\ldots,u_{r}),\ldots,(u_{2}+\alpha_{1}-\alpha_{2},u_{2}+1,u_{3}\ldots,u_{r})\}$$ belong to the same class in $V(\alpha_{1},\ldots,\alpha_{r})/\ker(\bar{f})$. 

As the cardinality of each class is $\alpha_1-\alpha_2+1$, one element is missing in this class.

 By (\ref{chain_un}), a class must be a chain. Therefore,   the last element $(u_{2}+\alpha_{1}-\alpha_{2},u_{2}+1,\ldots,u_{r})$ must have a father in the class. The possible fathers are among the following sequence of hypertuples: $(u_{2}+\alpha_{1}-\alpha_{2}+1,u_{2}+1,u_{3},\ldots,u_{r})$, $(u_{2}+\alpha_{1}-\alpha_{2},u_{2}+2,u_{3},\ldots,u_{r})$,\ldots, $(u_{2}+\alpha_{1}-\alpha_{2}+1,u_{2}+1,u_{3}+1,\ldots,u_{r}),\ldots$:
\begin{itemize}
  \item If $(u_{2}+\alpha_{1}-\alpha_{2},u_{2}+2,u_{3},\ldots,u_{r})\in [(u_{2}+1,u_{2}+1,u_{3},\ldots,u_{r}) ] $, then
    $$(u_{2}+1,u_{2}+1,u_{3},\ldots,u_{r})\vee(u_{2}+2,u_{2}+2,u_{3},\ldots,u_{i},\ldots,u_{r})\sim$$
    $$\sim (u_{2}+\alpha_{1}-\alpha_{2},u_{2}+2,u_{3},\ldots,u_{r})\vee (u_{2}+2,u_{2}+2,u_{3},\ldots,u_{i},\ldots,u_{r}),$$
   that is 
     $$ (u_{2}+2,u_{2}+2,u_{3},\ldots,u_{r})\sim (u_{2}+\alpha_{1}-\alpha_{2},u_{2}+2,u_{3},\ldots,u_{r}),$$
      but it means that $ (u_{2}+2,u_{2}+2,u_{3},\ldots,u_{r})\in [(u_{2}+1,u_{2}+1,u_{3},\ldots,u_{r})]$, and this is not possible because then the class will have more elements than  its cardinality.

    \item If $(u_{2}+\alpha_{1}-\alpha_{2},u_{2}+1,u_{3},\ldots,u_{i}+1\ldots,u_{r})\in [(u_{2}+1,u_{2}+1,u_{3},\ldots,u_{r})]$, then
    $$(u_{2}+1,u_{2}+1,u_{3},\ldots,u_{r})\vee (u_{2}+1,u_{2}+1,u_{3},\ldots,u_{i}+1,\ldots,u_{r})\sim$$
    $$\sim(u_{2}+\alpha_{1}-\alpha_{2},u_{2}+1,u_{3},\ldots,u_{i}+1,\ldots,u_{r})\vee (u_{2}+1,u_{2}+1,u_{3},\ldots,u_{i}+1,\ldots,u_{r}),$$
    that is 
     $$ (u_{2}+1,u_{2}+1,u_{3},\ldots,u_{i}+1,\ldots,u_{r})\sim (u_{2}+\alpha_{1}-\alpha_{2},u_{2}+1,u_{3},\ldots,u_{i}+1,\ldots,u_{r}),$$
     and this is not possible because the class will have more elements than its cardinality.
\end{itemize}
Therefore, the new element belonging to the class is $(u_{2}+\alpha_{1}-\alpha_{2}+1,u_{2}+1,u_{3},\ldots,u_{r})$.\\

Taking into account that $C_2, RC_2 \in V(\alpha_{1},\ldots,\alpha_{r})/\ker(\bar{f})$, by recurrence the result follows.

\end{proof} 
%%%%%%%%%%%%%%%%%%%%%%%%%%%%%%%%%%%%%%%%%%%%%%%%%%%%%%%%%%%%%%%%%%%%%%
%Lets start with an easy cases:
%\begin{remark}
 %Notice that
 %\begin{itemize}
%  \item If $r=s=1$, then
  %$$V(\alpha_{1})\simeq V(\beta_{1})\Leftrightarrow \alpha_{1}=\beta_{1}$$
  %\item If $r=1$ and $s=2$, then  $V(\alpha_{1})\simeq V(\beta_{1},\beta_{2})$ implies that all the hypertuples in $V(\beta_{1},\beta_{2})$ must to have an unique son an this implies that 
  %$\beta_{2}=\beta_{1}-1$, moreover $\alpha_{1}=\beta_{1}+\beta_{2}$ and therefore
%  $$V(\alpha_{1})\simeq V(\beta_{1},\beta_{2})\Leftrightarrow \beta_{2}=\beta_{1}-1,\alpha_{1}=2\beta_{1}-1$$
 % \end{itemize}
  %\end{remark}

\medskip

\medskip
 
These equivalence relations allow us to prove the following lemma.

\begin{proposition}\label{prop:quo} Let $V(\alpha_{1},\ldots,\alpha_{r})$ and $V(\beta_{1},\ldots,\beta_{r})$ be hyperlattices, then
 $$V(\alpha_{1},\ldots,\alpha_{r})\simeq V(\beta_{1},\ldots,\beta_{r}) \ \text{ if and only if } \
 \alpha_{i}=\beta_{i}, \ i=1, \ldots, r.$$
\end{proposition}
\begin{proof}
We apply induction. For $r=1$ and $r=2$ the result is proved in Propositions 
\ref{propo:r1} and \ref{propo:case22}.

For $r\geq3$, let $f:V(\alpha_{1},\ldots,\alpha_{r})\longrightarrow V(\beta_{1},\ldots,\beta_{r})$ be an isomorphism, then by Proposition~\ref{lemma:rCi},   $f(C_{1})=C'_{1}$. Hence, we analyze the following cases:
\begin{itemize}
 \item If $\alpha_{1}-\alpha_{2}>1$ and $\beta_{1}-\beta_{2}>1$, then $f(C_{2})=C'_{2}$, and  by  (\ref{eqganchito2}) and Theorem~\ref{teore:iso2}  we have 
 $$  V(\alpha_{1},\ldots,\alpha_{r})/\sim_{2}= V(\alpha_{1},\ldots,\alpha_{r})/\ker(\bar{f})\simeq V(\beta_{1},\ldots,\beta_{r})/\sim_{2}.$$
 Therefore, by Lemma~\ref{lemQuo} we obtain 
 $$V(\alpha_{2},\ldots,\alpha_{r})\simeq V(\beta_{2},\ldots,\beta_{r}).$$
 The result follows by induction hypothesis.
\item If $\alpha_{1}-\alpha_{2}=1$ and $\beta_{1}-\beta_{2}=1$, then $f(C_{3})=C_{3}'$, and by an analogous argument to that of the previous case we obtain 
$$V(\alpha_{3},\ldots,\alpha_{r})\simeq V(\beta_{3},\ldots,\beta_{r}).$$ The result follows by induction hypothesis.
 \item If $\alpha_{1}-\alpha_{2}>1$ and $\beta_{1}-\beta_{2}=1$, then  $f(C_{2})=C'_{3}$, which implies $\alpha_{1}-\alpha_{2}=2(\beta_{1}-\beta_{3})$ and $f(RC_{2})=RC'_{3}$.  It means  that $\alpha_{1}-\alpha_{2}+1=2(\beta_{1}-\beta_{2})$, and we have a contradiction.

\end{itemize}

\end{proof}

\subsection{Main theorem} \label{main}

For the general case the  result is summarized in the following main theorem, which is proved in the previous results. 

\begin{theorem}\label{teo:main}
 Let $A,B\in \mathbb{C}^{n\times n}$ be nilpotent matrices. Let  $\alpha=(\alpha_{1},\ldots,\alpha_{r})$ and $\beta=(\beta_{1},\ldots\beta_{s})$ be the reduced Segre characteristic  of $A$ and $B$, respectively. Then, $\Hinv(A)\simeq \Hinv(B)$ if and only if one of the following conditions are satisfied:
 \begin{enumerate}
  \item $\alpha=(5,2)$ and $\beta=(4,2,1)$ or vice-versa.
  \item $\alpha=(l,l-1)$ and $\beta=(2l-1)$ for some $l\geq 2$ or vice-versa.
  \item $\alpha=\beta$.
 \end{enumerate}
\end{theorem}
\section*{Acknowledgments}

The second author is partially supported by grant PID2019-104047GB-I00. The second and third authors are partially supported by grant PID2021-124827NB-I00 funded by MCIN/AEI/ 10.13039/501100011033 and by ``ERDF A way of making Europe'' by the ``European Union''.

\bibliographystyle{model1a-num-names}

\end{document}